\pdfoutput=1
\documentclass[12pt,a4paper,twoside]{scrartcl}
\usepackage[utf8]{inputenc}
\usepackage[english]{babel}
\usepackage{amsfonts}
\usepackage{amssymb}
\usepackage{amsthm}

\usepackage{mathtools}
\usepackage{dynkin-diagrams}
\usepackage{dsfont}
\usepackage{yfonts} % for German letters in math
\usepackage{upgreek} % upright greek letters
\usepackage{ytableau}
\usepackage{float}
\usepackage[left=3cm,right=3cm,top=3cm,bottom=3cm]{geometry}

\usepackage{imakeidx}
\makeindex[columns=2, title=Index, intoc]

\usepackage[refpage, intoc, english]{nomencl}
\makenomenclature

\usepackage{fancyhdr}
\usepackage{framed,color}
\usepackage{ifthen} %essential for logical dependencies in command definition
\usepackage{mathabx}
\usepackage{mathrsfs}
\usepackage{amsmath}
\usepackage[shortlabels]{enumitem} 

\usepackage{tikz}
\usetikzlibrary{arrows}
\usetikzlibrary[topaths]
\usepackage{cancel}

\usepackage{tikz-cd}
\usetikzlibrary{cd}
\usepackage{hyperref}
\usetikzlibrary{babel}
\usepackage[babel=true]{csquotes}

\usepackage{picture} %only so that \makebox(2cm,3cm){Text} stops making "errors" (despite working properly)

\usepackage{stackengine}
\usepackage{scalerel}

\usepackage{cleveref}
\makeatletter
\def\smallunderbrace#1{\mathop{\vtop{\m@th\ialign{##\crcr
				$\hfil\displaystyle{#1}\hfil$\crcr
				\noalign{\kern3\p@\nointerlineskip}%
				\tiny\upbracefill\crcr\noalign{\kern3\p@}}}}\limits}
\makeatother
\let\saveLongrightarrow\Longrightarrow
\makeatletter
\renewcommand*{\Longrightarrow}{%
	\mathrel{\rlap{\fontfamily{cmrx}\fontencoding{OT1}\selectfont=}%
		\hphantom{\saveLongrightarrow}%
		\llap{$\m@th\Rightarrow$}}}
\makeatother

\let\geq\geqslant
\let\leq\leqslant

\usepackage{amsthm}

\newcommand{\Res}{\mathrm{Res}}
\newcommand{\ct}{\mathrm{ct}}
\newcommand{\Ppl}{P^+(2\Sigma)}

\newcommand{\uq}{\mathbf {U}}
\newcommand{\uqcon}{\check{\mathbf {U}}}

\newcommand{\qfa}{R_q(G)}

\newcommand{\Q}{\mathbb{Q}}
\newcommand{\Z}{\mathbb{Z}}
\newcommand{\N}{\mathbb{N}}
\newcommand{\C}{\mathbb{C}}
\newcommand{\F}{\mathbb{F}}

\newcommand{\A}{\mathbf{A}}

\theoremstyle{definition}
\newtheorem{defi/}{Definition}[section]
\newenvironment{defi}
{%
	\pushQED{\qed}\begin{defi/}}
	{\popQED\end{defi/}}
\newtheorem{ex/}[defi/]{Example}
\newenvironment{ex}
{%
	\pushQED{\qed}\begin{ex/}}
	{\popQED\end{ex/}}
\newtheorem{notation/}[defi/]{Notation}
\newenvironment{notation}
{%
	\pushQED{\qed}\begin{notation/}}
	{\popQED\end{notation/}}
\theoremstyle{plain}
\newtheorem{prop}[defi/]{Proposition}
\newtheorem{thm}[defi/]{Theorem}

\newtheorem{cor}[defi/]{Corollary}
\newtheorem{lem}[defi/]{Lemma}

\theoremstyle{remark}
\newtheorem{rem}[defi/]{Remark} % instead of \newtheorem*
\usepackage{braket} %%% !!! this package is not in packages.tex !!!

\newcommand{\llbracket}{\left[\!\left[}
\newcommand{\rrbracket}{\right]\!\right]}

\renewcommand{\set}[1]{\left\lbrace#1\right\rbrace}

\makeatletter
\def\@tvsp{\mathchoice{{}\mkern-4.5mu}{{}\mkern-4.5mu}{{}\mkern-2.5mu}{}}
\def\ltrivert{\left|\@tvsp\left|\@tvsp\left|}
\def\rtrivert{\right|\@tvsp\right|\@tvsp\right|}
\def\lbivert{\left|\@tvsp\left|} 
\def\rbivert{\right|\@tvsp\right|} 
\def\llangle{\left\langle\@tvsp\left\langle}
\def\rrangle{\right\rangle\@tvsp\right\rangle}
\makeatother

\newcommand{\I}{\mathbf{I}}
\newcommand{\uqb}{\mathbf{B}_{\boldsymbol{c}}}
\newcommand{\uqbs}{\mathbf{B}_{\boldsymbol{c},\boldsymbol{s}}}
\newcommand{\uqds}{\mathbf{B}_{\boldsymbol{d},\boldsymbol{t}}}
\newcommand{\uqbl}{\mathbf{{B}}_{\boldsymbol{c},\boldsymbol{s}}^\omega}
\newcommand{\uqd}{\mathbf{B}_{\boldsymbol{d}}}

\newcommand{\uqbi}{\mathbf{B}_{\boldsymbol{c}_i,\boldsymbol{s}_i}}
\newcommand{\uqbis}{\mathbf{B}_{\boldsymbol{c}_i}}
\newcommand{\uqi}{\mathbf{U}_i}

\newcommand{\cl}{\mathsf{cl}}

\newcommand{\bottom}{\mathcal  B^+(\uq,\uqbs,\chi)}
\newcommand{\spqfa}{^{\uqbs^\rho}_{\quad\epsilon}\qfa_\epsilon^{\uqbs}}

\catcode`,\active

\catcode`\,12

\catcode`,\active

\catcode`\,12

\newcount\mycount

\newcommand{\Veranstaltung}{Algebraic Topology}
\newcommand{\Aufgabennr}{1}

\makeatletter %!!!
\fancypagestyle{plain}{% % <-- this is new
	\fancyhf{} 
	\fancyfoot[LE,RO]{\thepage} % same placement as with page style "fancy"
	}

\makeatother %!!!

\numberwithin{table}{section} % Table enumeration

\definecolor{shadecolor}{rgb}{0.8,0.8,0.8}

\newsavebox{\foobox}

\makeatletter
\newcommand{\thickbar}{\mathpalette\@thickbar}
\newcommand{\@thickbar}[2]{{#1\mkern1.5mu\vbox{
			\sbox\z@{$#1\mkern-.5mu#2\mkern-.5mu$}%
			\sbox\tw@{$#1\overline{#2}$}%
			\dimen@=\dimexpr\ht\tw@-\ht\z@-.8\p@\relax
			\hrule\@height.8\p@ % adjust for the desired rule thickness
			\vskip\dimen@
			\box\z@}\mkern1.5mu}
}
\makeatother

\newcommand{\myrule} [3] []{
	\begin{center}
		\begin{tikzpicture}
			\draw[#2-#3, ultra thick, #1] (0,0) to (0.65\linewidth,0);
		\end{tikzpicture}
	\end{center}
}

\newcommand{\demph}[2][]{\emph{#2}\ifthenelse{\equal{#1}{}}{\index{#2|hfill}}{\index{#1|hfill}}} %muss vllt überarbeitet werden

\makeatletter
\def\namedlabel#1#2{\begingroup
	#2%
	\def\@currentlabel{#2}%
	\phantomsection\label{#1}\endgroup
}
\makeatother

\usepackage[title]{appendix}

\usepackage{pdfpages}
\usepackage{comment}
\newcommand{\Addresses}{{% additional braces for segregating \footnotesize
		\bigskip
		\footnotesize
		\textsc{IMAPP, Radboud Universiteit, P.O. Box 9010, 6500 GL Nijmegen, The Netherlands}\par\nopagebreak
		\textit{E-mail address}: \texttt{stein.meereboer@ru.nl}
		
}}
\author{Stein Meereboer}
\title{Quantum spherical functions of type $\chi$ as Macdonald-Koornwinder polynomials}

\begin{document}

\makeatletter
\begin{center}
{\Huge \bfseries \textsf \@title}\\[6pt]
{\large\@author}\\
{\Addresses}

\end{center}
\makeatother
\noindent

\begin{abstract}
	The theory of quantum symmetric pairs is applied to $q$-special functions. Previous work shows the existence of a family $\chi$-spherical functions indexed by the integers for each Hermitian quantum symmetric pair. In this paper a distinguished family of such functions, invariant under the action of the Weyl group of the restricted roots, is shown to be a family of Macdonald-Koornwinder polynomials if the restricted root system is reduced or if the Satake diagram is of type $\mathsf{AIII_a}$.
\end{abstract}
%\tableofcontents
\numberwithin{equation}{section}
\section{Introduction and overview}
A beautiful result of Gail Letzter shows that quantum zonal spherical functions can be realized as Macdonald polynomials \cite{Letzter2004}. Her theory and results unify many case-by-case computations of Noumi, Sugatani, Dijkhuizen, Stokman, Koornwinder and others \cite{Noumi1992} \cite{Noumi1995}, \cite{Noumi1996}, \cite{Koornwinder1993}, \cite{Stok98}. The connection between zonal spherical functions and $q$-special functions has examples in the more general setting of quantum spherical functions of type $\chi$. Koelink \cite[\mbox{Thm 7.1}]{Koelink1996} and Noumi and Mimachi \cite[\mbox{Thm 4}]{Noumi1992} show that quantum spherical functions of type $\chi$ on the quantum analog of the symmetric pair $(SU(2),U(1))$ can be interpreted as Askey-Wilson polynomials. 

Classically, spherical functions of type $\chi$ can be interpreted as hypergeometric functions, cf. \cite[\mbox{Thm 5.2.2}]{Heckman1995}. This motivates the question whether quantum spherical functions of type $\chi$ can be realized as $q$-hypergeometric functions. The main result of this paper gives as an affirmative answer to this question for all quantum symmetric pairs with reduced restricted root system and in type $\mathsf{AIII_a}$. The approach avoids any radial part calculations. Instead we rely on the identification in the zonal spherical case, which allows us to explicitly determine the orthogonality measure for quantum spherical functions of type $\chi$. Our identification avoids extensive case work. 

As a consequence, we provide simple, uniform, formulas for the parameters appearing in the Macdonald-Koornwinder polynomials. The identification allows us to show properties of the orthogonal polynomials. With use of the $\imath$-bar involution of Bao and Wang \cite{Bao2018} and the relative braid group action of Wang and Zhang \cite{Wang2023}, we recover the known $q\to q^{-1}$ symmetry of Koorwinder-Macdonald polynomials. Furthermore, in the reduced case we identify the Macdonald $q$-difference operators as operators arising form the center of the quantum group. Finally, tensor product decompositions yield relations for connection coefficients of Macdonald-Koornwinder polynomials.
\\
\\
This paper is organized as follows.\\
In Section \ref{sec:QG} we give preliminaries about quantum groups, the quantized function algebra and quantum symmetric pairs including a review of character theory for Hermitian quantum symmetric pairs, cf. \cite{Meereboer2025}.

Section \ref{sec:OPSF} introduces our primary objects of study: spherical functions of type $\chi$, and Macdonald-Koornwinder polynomials. We also discuss the identification of zonal spherical functions with orthogonal polynomials, cf. \cite{Letzter2004}, \cite{Noumi1996}.

In Section \ref{sec:Const} we construct, for each family of spherical functions of type $\chi$, a distinguished family of orthogonal polynomials. The essence of this construction is as follows: by viewing the space of $\chi$-spherical functions as a free module of rank one over the algebra of zonal spherical functions, we associate zonal spherical functions to each family of $\chi$-spherical functions. Using the Haar integral on the quantized function algebra and quantum Schur orthogonality we show that the polynomials are orthogonal with respect to a weight involving the orthogonality weight for zonal spherical functions and the fundamental spherical function of type $\chi$. Here the fundamental spherical functions are certain minimal elements in the space of spherical functions. Determining these polynomials thus reduces to identifying the fundamental spherical function.

Section \ref{sec:fund} determines fundamental spherical functions by reducing to the rank one cases $\mathsf{AI}$ and $\mathsf{AIV}$. 

Combining the results of Sections \ref{sec:Const} and \ref{sec:fund}, Section \ref{sec:iden} identifies the constructed orthogonal polynomials with Macdonald-Koornwinder polynomials based on the explicit expression of the weight for which the polynomials are orthogonal. 

Finally, Section \ref{sec:apl} applies the identification with orthogonal polynomials and quantum group techniques to establish three results for the orthogonal polynomials. First, we prove the invariance of orthogonal polynomials under sending $q\to q^{-1}$. Second, the Macdonald-Koornwinder $q$-difference operators are identified with central elements of the quantum group. Third, we derive relations for the connection coefficients of Macdonald-Koornwinder polynomials.
\begin{center}
	$\textsc{Acknowledgment}$
\end{center}
The author likes to thank Erik Koelink for his guidance and for many valuable comments. The author likes to thank Maarten van Pruijssen for the discussions and many valuable comments. The author also likes to thank Philip Schl\"osser for the discussions regarding Macdonald-Koornwinder polynomials, pointing out typos in an earlier version, and, in particular, for the details of the proof of Lemma  \ref{lem:constant}. The author is grateful to the
referees for a careful reading, suggestions and corrections.
\begin{center}
	$\textsc{Funding}$
\end{center}
My research is funded by NWO grant \texttt{OCENW.M20.108}. 
\section{Quantum groups and quantum symmetric pairs}\label{sec:QG}
We introduce notation as in \cite{Lusztig2010}. Let $\mathfrak{g}$ be be a complex semisimple Lie algebra with  Cartan matrix $C=(c_{ij})_{i,j\in \I}$. Let $D=\text{diag}(d_i\,:\,d_i\in \Z_{\geq 1}, \,i\in\I)$ be a symmetrizer, meaning that $DC$ is symmetric and $\gcd\{d_i\,:\, i\in\I\}=1$. We fix a set of simple roots $\Pi=\{\alpha_i\,:\,i\in \I\}$ and a set of simple coroots $\Pi^\vee=\{\alpha^\vee_i\,:\,i\in \I\}$ with corresponding root and coroot systems $\mathcal R$ and $\mathcal R^\vee$. The root lattice is denoted by $\Z\I=\oplus_{i\in \I } \Z\alpha_i$. The root lattice is equipped with the normalized Killing form $(\,,\,)$, scaled so that short roots have length $2$. The Weyl group $W$ is generated by the simple reflections $s_i:\Z\I\to \Z\I$, for $i\in \I$ defined by $s_i(\alpha_j)=\alpha_j-c_{ij}\alpha_i,$ for $j\in \I$. 
\\

Let $q$ be an indeterminate, and let $\Q(q)$ denote the field of rational functions in $q$ with coefficients in $\Q$. Denote by $\F$ the algebraic closure of $\Q(q)$. We fix a group
homomorphism $\Q \to \F^\times$ mapping $1 \to q$. This corresponds to a consistent choice of rational
powers $q^x$ of $q$ $(x \in \Q)$. The field $\F$ is equipped with a unique $\Q$-linear automorphism $\overline{\,\cdot \,}:\F\to \F$ such that $q$ is mapped to $q^{-1}$. Let $\A$ be the smallest unital subring of $\F$ for which $q$ is contained in $\A$ and such that every element in $\A$ that is not contained in the ideal generated by $(q - 1)$ has a square root and is invertible, cf. \cite[\mbox{\textsection 1.}]{Letzter2000}. Define $\cl:\A\to \A/(1-q)\A\cong \C$ as the map $f\mapsto f(1)$. For each $i\in\I$ and $n\in \Z_{\geq0}$, set
$$q_i:= q^{d_i},\qquad [n]_i:=\cfrac{q_i^n-q_i^{-n}}{q_i-q_i^{-1}}\qquad\text{and}\qquad [n]_{i}!:=\prod_{k=1}^n[k]_i.$$

Let $(Y,X, \langle \,,\, \rangle.\dots)$ be simply connected root datum of type $(\I,\cdot)$, cf. \cite[\mbox{2.2}]{Lusztig2010}.  By definition, there are embeddings $\I\to X$, $i\mapsto \alpha_i$ and $\I\to Y$, $i\mapsto \alpha_i^{\vee}$, such that $Y=\Z [\alpha_i^\vee]$ and $X=\mathrm{Hom}_\Z(Y,\Z)$, respectively. The pairing $\langle\,,\,\rangle$ and the form $(\,,\,)$ are related by the formula $\langle \alpha_i^\vee,\alpha_j\rangle=\frac{2(\alpha_i,\alpha_j)}{(\alpha_i,\alpha_i)}$, for $i,j\in \I$. The lattice $X$ is partially ordered by the rule
$$\lambda\leq\lambda'\quad\text{if and only if}\quad \lambda'-\lambda\in \Z_{\geq0}\I.$$

The quantum group $\uq$, associated to the root datum, is the unital associative algebra over $\F$, generated by the symbols $E_i,F_i$ for $i\in \I$ and $K_{h}$ for $h\in Y$. These generators satisfy the following relations for $ i,j\in\I,$ and $ h,h_1,h_2\in Y:$
\begin{align*}
	K_0&=1,\\
	K_{h_1}K_{h_2}&=K_{h_1+h_2},\\
	K_{h}E_i&= q^{\langle h,\alpha_i\rangle}E_iK_h,\\
	K_{h}F_i&= q^{-\langle h,\alpha_i\rangle}F_iK_h,\\
	E_iF_j-F_jE_i&=\delta_{i,j}\cfrac{K_i-K_i^{-1}}{q_i-q_i^{-1}},
\end{align*}
as well as the quantum Serre relations for $i\neq j$,
\begin{align*}
	\sum_{r+s=1-a_{i,j}}(-1)^sE_i^{(r)}E_jE_i^{(s)}=\sum_{r+s=1-a_{i,j}}(-1)^s F_i^{(r)}F_jF_i^{(s)}=0.
\end{align*}
Here we use the notation
$$K_i:=K_{d_i\alpha_i^\vee},\qquad E_i^{(n)}:=\frac{1}{[n]_i!}E_i^n\qquad\text{and}\qquad F_i^{(n)}=\frac{1}{[n]_i!}F_i^n.\qedhere$$
The quantum group $\uq$ has the structure of a Hopf-algebra $(\uq, \triangle,\epsilon,\iota,S)$, cf. \cite[\mbox{Ch 4}]{Jantzen1996}. The conventions for the antipode, coproduct and counit are
    \begin{align*}
	\triangle(E_i)&= E_i\otimes 1 + K_i\otimes E_i&	\triangle(F_i)&= F_i\otimes K_i^{-1}+ 1\otimes F_i&	\triangle(K_h) &= K_h\otimes K_h\\
	S(E_i) &= -K_i^{-1}E_i	&S(F_i) &= -F_iK_i	&S(K_h) &= K_{-h}\\
	\epsilon(E_i)&=\epsilon(F_i)= 0&	\epsilon(K_h)&= 1
\end{align*}
for $i\in I$ and $h\in Y$.

We use \textit{Sweedler notation} to denote the coproduct, writing
\begin{equation}\label{eq:sweed}\triangle(x)=x_{(1)}\otimes x_{(2)}\in \uq\otimes \uq.\end{equation}
\begin{enumerate}[(i)]
	\item Let $\omega: \uq\to \uq$ denote the unique algebra automorphism that satisfies 
	\begin{equation}\label{eq:chevaut}
		E_i\mapsto F_i,\qquad F_i\mapsto E_i,\qquad K_i\mapsto K_i^{-1},\qquad \text{ for all }i\in \I.
	\end{equation}
	\item Let $\varrho: \uq\to \uq$ denote the unique algebra anti-automorphism that satisfies 
	\begin{equation}\label{eq:rhoaut}
		E_i\mapsto F_i,\qquad F_i\mapsto E_i,\qquad K_i\mapsto K_i,\qquad \text{ for all }i\in \I.
	\end{equation}
	\item Let $\overline{\,\cdot\,}:\uq\to \uq$ denote the unique $\Q$-algebra automorphism that satisfies
	\begin{equation}\label{eq:barinv}
		E_i\mapsto E_i,\qquad F_i\mapsto F_i,\qquad K_i\mapsto K_i^{-1},\qquad q\mapsto q^{-1}\qquad \text{ for all }i\in \I
	\end{equation}
	\item For each $\boldsymbol{a}\in (\F^\times)^\I$ let $\Phi_{\boldsymbol{a}}:\uq\to\uq$ denote the Hopf-algebra automorphism that satisfies
	\begin{equation}\label{eq:aut}
		E_i\mapsto a_i^{1/2}E_i,\qquad F_i\mapsto a_i^{-1/2}F_i,\qquad K_i\mapsto K_i,\qquad\text{ for all }i\in \I.
	\end{equation}
\end{enumerate}

For each $i\in \I$, let $T_i:\uq\to \uq$ denote the algebra automorphism corresponding to $T_{i,+1}''$ as defined in \cite[\mbox{37.1}]{Lusztig2010}. 
The operators $T_i$ satisfy the braid relations. That is to say, if $w = s_{i_1} \dots s_{i_n} \in W$ is a reduced expression, there exist a well-defined automorphism $T_{w}=T_{s_{i_1}} \dots T_{s_{i_n}} $ of $\uq$.
\subsection{Finite dimensional weight modules}\label{subsec:intromod}
By convention $\uq$-modules are left modules, unless stated otherwise.
A $\uq$-module $M$ is called a \textit{weight $\uq$-module} if $M$ has a decomposition of the form
$$M=\bigoplus _{\lambda\in X}M_{\lambda },\qquad\text{where}\qquad M_\lambda=\{m\in M \,:\, K_h m= q^{\langle h,\lambda \rangle}m, \,\text{for all } h\in Y\}.$$
We say that a vector $v\in M$ is a \textit{weight vector of weight $\lambda$} if $K_h v= q^{\langle h,\lambda \rangle}v$ for each $ h\in Y.$
A weight $\uq$-module $M$ is called a \textit{highest weight module} if there exists a weight vector $v\in M$ such that $\uq v=M$ and $E_iv=0$ for all $i\in \I.$ 
The irreducible finite dimensional weight $\uq$-modules are classified by the dominant weights $X^+$, cf. \cite[\mbox{3.5}]{Lusztig2010}. 
The simple $\uq$-module corresponding to a dominant weight $\lambda\in X^+$ is denoted by $L(\lambda)$. 
The theory of highest weight $\uq$-modules has a natural analog for right $\uq$-modules, cf. \cite[\mbox{3.5.7}]{Lusztig2010}. 

Let $L(\lambda)$ be an irreducible weight $\uq$-module and $v_\lambda$ be a highest weight vector. Then, there exists a \textit{bar involution} $\overline{\,\cdot\,}^{L(\lambda)}: L(\lambda)\to L(\lambda)$ relative to $v_\lambda$ satisfying
$$\overline{Xv_\lambda}^{L(\lambda)}:=\overline{X} v_\lambda\qquad\text{for}\qquad X\in \uq.$$
This is well-defined by \cite[\S19.3.4]{Lusztig2010}.
If it is clear from the context, we write $\overline{\,\cdot\,}^{L(\lambda)}=\overline{\,\cdot\,}$. Analogously, irreducible right weight $\uq$-modules are equipped with bar involutions relative to chosen lowest weight vector.

If $\vartheta:\uq\to \uq$ is an (anti)-automorphism and $M$ is a left $\uq$-module, we define $^\vartheta M$ to be $\uq$-module with underlying vector space $M$ and such that the action is twisted by $\vartheta$. 
%To be precise, the action of $\uq$ on $^\vartheta M$ is defined by
%\[x \triangleright m=\vartheta(x)m,\qquad \text{for }x\in \uq, m\in M.\]
Thus, if $\vartheta$ is an automorphism $^\vartheta M$ is a left $\uq$-module, and if $\vartheta$ is an anti-automorphism $^\vartheta M$ is a right $\uq$-module.
\subsection{The quantized function algebra and bilinear forms}\label{subsec:qfa}
The \textit{quantized function algebra} $\qfa$ is the subspace of $\uq^{\ast}=\mathrm{Hom}(\uq,\F)$ spanned by matrix coefficients for finite dimensional weight $\uq$-modules.
The quantized function algebra has the structure of a
Hopf-algebra $(\qfa, \triangle,\epsilon,\iota,S)$ which yields a Hopf-algebra duality between $\uq$ and $\qfa$, cf. \cite[\mbox{Ch 7}]{Jantzen1996}. 

For each dominant weight $\lambda\in X^+$, the dual $L(\lambda)^\ast$ becomes a right $\uq$-module by defining
$$(fX)(v):=f(Xv)\qquad \text{where}\qquad X\in \uq,\, f\in L(\lambda)^\ast,\, v\in L(\lambda).\qedhere$$
Let $f\otimes v\in L(\lambda)^\ast\otimes L(\lambda)$, then we write the corresponding element of $\qfa$ as $c^{L(\lambda)}_{f,v}.$ If the module is clear from the context, we abbreviate this as $c_{f,v}$. 
%w∗ ,w
%We define $c^{L(\lambda)}:L(\lambda)^\ast\otimes L(\lambda)\to \qfa$ by
%$$c^{L(\lambda)}_{f,v}(X):=f(X v)\qquad\text{where}\qquad f\in L(\lambda)^\ast,v\in L(\lambda), \, X\in \uq.$$
The quantized function algebra $\qfa$ has the quantum Peter-Weyl decomposition
$$\qfa=\bigoplus_{\lambda\in X^+} L(\lambda)^\ast\otimes L(\lambda)$$
as a $\uq$-bimodule, cf. \cite[\mbox{7.2}]{Kashiwara1993}. For each weight $\lambda\in X^+$ the $\uq$-bimodule 
\[L(\lambda)^\ast\otimes L(\lambda)\]
is referred to as a \textit{Peter-Weyl block.}

Next we describe an explicit realization of $L(\lambda)^\ast$ via the Shapovalov form. 
Define $L^r(\lambda)$ as the simple twisted right module $^\varrho L(\lambda)$, where $\varrho$ is the anti-automorphism defined in (\ref{eq:rhoaut}).
As observed by Masaki Kashiwara in \cite[\mbox{7.1.3}]{Kashiwara1993}, for each dominant weight $\lambda\in X^+$, the module $L(\lambda)$ is equipped with a unique non-degenerate symmetric bilinear form, relative to a chosen highest weight vector $v_\lambda \in L(\lambda)$, satisfying $(v_{\lambda},v_{\lambda})=1$ and 
\begin{equation}\label{eq:Shapovalof}
	(v,Xw)=(\varrho(X)v,w)\qquad \text{where}\qquad v,w\in L(\lambda),\, X\in\uq,
\end{equation}
where $\varrho$ is the anti-automorphism defined in \eqref{eq:rhoaut}.
This bilinear form, known as the \textit{Shapovalov form}, provides an explicit isomorphism $L(\lambda)^\ast\cong L^r(\lambda)$ as right $\uq$-modules, cf. \cite[\mbox{Prop 7.2.2}]{Kashiwara1993}.
\subsection{Quantum symmetric pairs}\label{sec:QSP}
In this section we introduce quantum symmetric pair coideal subalgebras following \cite{Letzter1999} and \cite{Kolb2014}. 

For each subset $\I_\bullet\subset \I$, let $w_\bullet$ denote the longest element in the parabolic subgroup $W_\bullet=\langle s_i:i\in \I_\bullet\rangle\subset W$.
\begin{defi}[Satake diagram]
A \textit{Satake diagram} $(\I_\bullet,\tau)$ consists of a subset $\I_\bullet\subset \I$ and an involutive diagram automorphism $\tau:\I\to \I$ with
\begin{enumerate}[(i)]
	\item $\tau|_{\I_\bullet}=-w_\bullet$;
	\item If $i\in \I\setminus\I_\bullet$ and $\tau(i)=i$ then $\alpha_i(\rho_\bullet^\vee)\in \Z$.
\end{enumerate}
Here $\rho_\bullet^\vee$ denotes the half sum of positive coroots relative to $\I_\bullet$. 
\end{defi}
Given a Satake diagram $(\I,\tau)$, we denote $\I_\circ:=\I\setminus\I_\bullet$ and set $$\I_{\mathrm{ns}}:=\{i\in \I_\circ\,:\, \tau(i)=i,\,\langle \alpha_i^\vee,\alpha_j\rangle =0\text{ for all }j\in \I_\bullet\}.$$
Moreover, we associate an involution $\Theta=-w_\bullet\circ\tau$ acting on $ X,Y$ and $\Z\I.$
Denote $\uq_\bullet= \F\langle E_j, F_j, K_j^\pm\,:\, j\in \I_\bullet\rangle$ the quantum group assosiated to $\I_\bullet$ and denote $\uq_\Theta^0=\F\langle K_h\,:\, h\in Y,\, \Theta(h)=h\rangle $ the $\Theta$-invariants of $\uq^0=\langle K_h:\,h\in Y\rangle $.
\begin{defi}[Quantum symmetric pair coideal subalgebra]
	Let $\uqbs$ denote the \textit{quantum symmetric pair coideal subalgebra} associated to the Satake diagram $(\I_\bullet,\tau)$. It is the subalgebra of $\uq$ generated by $\uq_\bullet$, $\uq_\Theta^0$, and the elements
	$$B_i^{\boldsymbol{c},\boldsymbol{s}}=B_i= F_i+c_i T_{w_\bullet}(E_{\tau(i)})K_i^{-1}+s_iK_i^{-1},\qquad\text{for all}\qquad i\in \I_\circ.$$
	We refer to $(\uq,\uqbs)$ as a quantum symmetric pair (QSP), and to the algebra $\uqbs$ as a QSP-coideal subalgebra of $\uq$.
	The QSP-coideal subalgebra $\uqbs$ depends on a family of scalars $(\boldsymbol{c},\boldsymbol{s})\in (\F^\times)^{\I_\circ}\times \F^{\I_\circ}$ that satisfy 
	{\belowdisplayskip=-12pt\begin{align*}
		\boldsymbol{c}\in \mathcal C&= \{\boldsymbol{d}\in(\F^\times)^{\I_\circ}: d_i=d_{\tau(i)}\quad \text{if}\quad \langle \alpha_i^\vee,\Theta(\alpha_i)\rangle=0  \}\\
		\boldsymbol{s}\in \mathcal S&= \{\boldsymbol{t}\in \F^{\I_\circ}: t_j\neq 0\implies (j\in \I_{\mathrm{ns}}\text{ and  }c_{ij}\in -2\Z_{\geq0}\,\forall i\in \I_{\mathrm{ns}}\setminus \{j\})\}.\hfill
	\end{align*}\qedhere}
\end{defi}
By \cite[Prop 5.2]{Kolb2014} QSP-coideal subalgebras $\uqbs$ are right coideal subalgebras, that is to say
$$\triangle (\uqbs)\subset \uqbs\otimes \uq.$$
For $\boldsymbol{s}=0$,  the corresponding QSP-coideal subalgebra is denoted by $\uqb$ and is said to be \textit{standard.}

For each tuple $\boldsymbol{a}\in (\F^\times)^\I,$ the coideal subalgebra $\Phi_{\boldsymbol{a}}(\uqbs)$ is a QSP-coideal subalgebra, here $\Phi_{\boldsymbol{a}}$ is the automorphism from (\ref{eq:aut}). 
Thus, the group $(\F^\times)^{\I}$ acts on the set of parameters $(\mathcal C,\mathcal S)$. 

A parameter $(\boldsymbol{c},\boldsymbol{s})$ is called \textit{specializable} if $(\boldsymbol{c},\boldsymbol{s})\in \A^{\I_\circ}\times \A^{\I_\circ}$ and  $\frac{c_i}{c_{\tau(i)}}+\A=(-1)^{\langle 2\rho_\bullet^\vee,\alpha_i\rangle }+\A$ for each $i\in \I_\circ$.
And a parameter $(\boldsymbol{c},\boldsymbol{s}) $ is called \textit{quasi specializable} if it is in the $(\F^\times)^{\I}$-orbit of a specializable parameter and $(\boldsymbol{c},\boldsymbol{s})\in \Q(q)^{\I_\circ}\times \Q(q)^{\I_\circ}$.

A Satake diagram $(\I_\bullet,\tau)$ is said to be \textit{irreducible} if for each $i,j \in \I$, there exists
a sequence $i= i_1\dots ,i_r = j \in \I$ such that for each $k = 1,...,r- 1$, we have either
$a_{i_k ,i_{k+1}} \neq 0$ or $i_{k+1} = \tau(i_k)$.

\begin{rem}\label{rem:asum}
Throughout this paper, we assume that our parameters are quasi-specializable and that our Satake diagrams $(\I_\bullet,\tau)$ are irreducible.
\end{rem}

\begin{defi}[Shift of base point]\label{def:shift}
For each character $\chi \in \widehat{\uqbs}$ with $\chi|_{\uq_\bullet}=\epsilon$ we define the algebra isomorphism
\begin{equation}\label{eq:shiftof}
	\Psi_\chi: \uqbs\to \uqds,\qquad b\mapsto \chi(b_{(1)})b_{(2)}\qquad b\in \uqbs.
\end{equation}
It is referred to as a \textit{shift of base point}, cf. \cite[\mbox{(10.1)}]{Kolb2023}. The isomorphism $\Psi_\chi$ maps $\uqbs$ to $\Psi_\chi(\uqbs)$, where $\Psi_\chi(\uqbs)=\uqds$ is the QSP-coideal subalgebra with parameters
\[d_i= \chi(K_i^{-1}K_{\tau(i)})c_i,\qquad t_i=\chi(B_i^{\boldsymbol{c},\boldsymbol{s}}),\qquad \text{where}\qquad i\in \I_\circ.\qedhere\]
\end{defi}

%\begin{rem}
%The reader is cautioned to distinguish the half sum of the positive roots $\rho$ and the shift of base point $\Psi_\chi$.
%\end{rem}
%The Satake diagrams $(\I_\bullet,\tau)$ are in bijection with the Satake diagrams that arise from classification of real simple Lie algebras. 

For each root $\alpha\in \mathcal R$, set $\widetilde{\alpha}= \frac{\alpha-\Theta(\alpha)}{2}$. The set 
\begin{equation}\label{eq:restricted}
	\Sigma=\{\widetilde{\alpha}\,:\, \alpha\in \mathcal R,\,\Theta(\alpha)\neq \alpha\}\end{equation}
is called the \textit{restricted root system}. The corresponding Weyl group is called the \textit{relative Weyl group} and has an embedding in $W$, cf. \cite[\mbox{Sec 2.3}]{Dobson2019}.

\begin{defi}\label{def:rank}
The number of $\tau$-orbits in $\I_\circ$ is called the \textit{rank} of the Satake diagram.
\end{defi}
Essential to many of our arguments is the reduction to the rank one case.
\begin{defi}[Rank one reductions]\label{def:rankone}
To each $i\in I_\circ$ we associate two lattices
\[Y_i=\Z[\alpha^\vee_j: j\in {\I_\bullet\cup \{i,\tau(i)\}}]\quad\text{and}\quad X_i:=\mathrm{Hom}_\Z(Y_i,\Z).\]
We write $\uqi$ for the quantized enveloping algebra associated to the simply connected root datum $(Y_i,X_i,\langle\,,\,\rangle,\dots)$. We set \[\Res(Y_i): X=\mathrm{Hom}_\Z(Y,\Z)\to \mathrm{Hom}_\Z(Y_i,\Z)=X_i,\qquad \lambda\mapsto \lambda|_{Y_i}\]
equal \textit{restriction of weights} from $X$ to $Y_i$.

To each $i\in I_\circ$ we associate the datum of a \textit{rank one Satake diagram} $(\I_\bullet,\tau|_{\I_\bullet\cup \{i,\tau(i)\}})$, by restricting, if necessary, to the connected component of the Dynkin diagram $\I_\bullet\cup \{i,\tau(i)\}$ containing $\{i,\tau(i)\}$. Let $\Sigma_i\subset \Sigma$ denote the restricted root system associated to the rank one Satake diagram. 

We denote $\uqbi\subset \uqi$ for the associated \textit{rank one QSP-coideal subalgebra}. Both $\uq_i$ and $\uqbi$ have a natural embedding
\[\uq_i\subset \uq\qquad \text{and}\qquad \uqbi\subset \uqbs.\qedhere\] 
\end{defi}
\begin{ex}\label{ex:1}
	\begin{enumerate}[(i)]
		\item The rank one reductions for the Satake diagram $\dynkin [labels*={1,2,3,4,5}]B{o2***}$ of type $\mathsf{BI}_5$ are the diagrams $\dynkin [labels*={1}]A{o}$ and $\dynkin [labels*={2,3,4,5}]B{o***}$.
		\item The rank one reductions for the Satake diagram $\dynkin [labels*={1,2,3,4,5},  involutions={[out=-20,in=20] {4}{5}}] D{*o*oo}$ of type $\mathsf{DIII_b}$ are the diagrams
		$\dynkin[labels*={1,2,3}] {A}{*o*}$, and $\dynkin[labels*={4,3,5}, involutions={[out=-60,in=-120] {1}{3}}] {A}{*o*}$. 
		\item For future reference, we refer to the rank one diagram $\dynkin[labels*={1,,,n}, involutions={[out=-60,in=-120] {1}{4}}] {A}{o*.*o}$ as being of type $\mathsf{AIV_n}$. For example, the rank two Satake diagram $\mathsf{DIII_b}$ has a \textit{rank one subdiagram} of type $\mathsf{AIV_3}$.\qedhere
	\end{enumerate}
\end{ex}
\subsection{Integrable characters of $\uqbs$}
In this section we review the theory of integrable characters of $\uqbs$. By a character we mean a unital $\F$-algebra homomorphism $\chi:\uqbs\to \F$. We write $\widehat{\uqbs}$ for the set of characters $\chi:\uqbs\to \F$
\begin{notation}
	Let $\chi:\uqbs\to \F$ be a character. Then we denote $V_{\chi}=\mathrm{span}_{\F}\{1\}$ the $\uqbs$-module with
	\[b1=\chi(b)1,\qquad \text{where}\qquad b\in \uqbs\qedhere\]
\end{notation}
\begin{defi}[Integrable character]
A character $\chi\in\widehat{\uqbs}$ is called \textit{integrable} if $\dim\big(\mathrm{Hom}_{\uqbs}( V_\chi,L(\lambda))\big):=[L(\lambda)|_{\uqbs}: V_\chi]\neq 0$ for some $\lambda\in X^+$. If $[L(\lambda)|_{\uqbs}: V_\chi]\neq 0$, we say that $L(\lambda)$ contains a module of type $\chi$.
\end{defi}
Let $i\in\I_\circ$ and recall from Definition \ref{def:rankone} the associated rank one Satake diagram $(\I_\bullet,\tau|_{\I_\bullet\cup \{i,\tau(i)\}})$.
\begin{defi}[Hermitian type \& Marked indices]  \cite[\S4]{Meereboer2025} \label{def:hermitian}
A QSP-coideal subalgebra $\uqbs$ is said to be of \textit{Hermitian type} if its corresponding Satake diagram is one of the Table \ref{table:1}
{\allowdisplaybreaks\begin{table}[H]
		\centering
		\begin{tabular}{|c | c | c |}
			\hline
			Type	& Satake diagram & Restricted root system   \\
			\hline
			$\mathsf{AIII_a}$	&$\begin{dynkinDiagram} A{IIIa}\dynkinRootMark{t}3 \dynkinRootMark{t}8
			\end{dynkinDiagram}$  & $\mathsf{BC_n}$   \\
			\hline
			$\mathsf{AIII_b}$	& $\begin{dynkinDiagram} A{IIIb} \dynkinRootMark{t}4
			\end{dynkinDiagram}$ & $\mathsf{C_n}$    \\
			\hline
			$\mathsf{BI}	$&  $\begin{dynkinDiagram} B{o2*.**}\dynkinRootMark{t}1\end{dynkinDiagram}$& $\mathsf{B_2}$   \\
			\hline
			$\mathsf{CI}$	& $\begin{dynkinDiagram} CI \dynkinRootMark{t}4\end{dynkinDiagram}$ &  $\mathsf{C_n}$, $n\geq2$ \\
			\hline
			$\mathsf {DI}$	& $\begin{dynkinDiagram} D{o2*.***}\dynkinRootMark{t}1\end{dynkinDiagram}$ &  $\mathsf{B_2}$  \\
			\hline
			$\mathsf{DIII_b}$	& $\begin{dynkinDiagram} D{IIIb} \dynkinRootMark{t}6\dynkinRootMark{t}7\end{dynkinDiagram}$ & $\mathsf{BC_n}$, $n\geq 2$   \\
			\hline		
			$\mathsf{EIII}$	&  $\begin{dynkinDiagram} E{III}\dynkinRootMark{t}1\dynkinRootMark{t}6\end{dynkinDiagram}$& $\mathsf{BC_2}$   \\
			\hline
			$\mathsf{EVII}$	&  $\begin{dynkinDiagram} E{VII} \dynkinRootMark{t}7\end{dynkinDiagram}$&  $\mathsf{C_3}$  \\
			\hline
		\end{tabular}
		\caption{Satake diagrams of Hermitian type} \label{table:1}
\end{table}}
In Table \ref{table:1} the mark $\scriptscriptstyle\otimes$ indicates an $i\in\I_\circ$ such that the rank one Satake diagram $(\I_\bullet,\tau|_{\I_\bullet\cup \{i,\tau(i)\}})$ is either of type $\mathsf{AIII_a}$ or that $i\in I_{\mathrm{ns}}$ and $a_{ij}\in 2\Z$ for all $j\in I_{\mathrm{ns}}$.
We denote $\I_{\scriptscriptstyle\otimes}\subset \I_\circ$ the set of \textit{marked indices}. 
\end{defi}
We note that for each irreducible type in Table \ref{table:1} there exists precisely one $\tau$-orbit in $\I_{\scriptscriptstyle\otimes}$.
We recall that the root system $\mathsf{B_2}$ is isomorphic to $\mathsf{C_2}$. Hence each restricted root system in Table \ref{table:1} is of type of type $\mathsf{C}$ or $\mathsf{BC}$. We borrow the following notation from \cite[\mbox{5.1.3.}]{Heckman1995}.
\begin{rem}\label{rem:ell}
For type $\mathsf{BC}$ we set $\Sigma_{\ell}$ equal to the long roots, $\Sigma_{m}$ equal to the middle short roots and $\Sigma_s=\frac{1}{2}\Sigma_{\ell}$. For type $\mathsf{C}$ we set $\Sigma_{\ell}$ equal to the long roots and $\Sigma_m$ equal to the short roots. This is also explained in \eqref{eq:bc}.
%Furthermore, in either case write $\Sigma_m$ for the middle long roots; these correspond to the short roots in $\Sigma$ in type $\mathsf{C}$ and correspond to the middle long roots in $\Sigma$ in type $\mathsf{BC}$.
%This notation is used in \cite[\mbox{5.1.3.}]{Heckman1995}. Consider the case that $\Sigma$ is of type $\mathsf{C}$. Then, $\Sigma_s$ is not contained in $\Sigma$.
\end{rem}
Recall from Remark \ref{rem:asum} that the Satake diagram $(\I_\bullet,\tau)$ is assumed to be irreducible. This implies that the associated symmetric pair $(\mathfrak{g},\mathfrak{g}^\Theta)$ is irreducible, i.e that $\mathfrak{g}$ cannot be written as the
direct sum of two semisimple Lie subalgebras which both admit $\Theta$ as an involution. Table \ref{table:1} is in one-to-one correspondence with the irreducible symmetric pairs $(\mathfrak{g},\mathfrak{g}^\Theta)$, cf. \cite[\mbox{X. \S6.3}]{Helga98}.

We can also characterize being Hermitian through the existence of nontrivial integrable characters.
\begin{prop}\cite[\mbox{Prop 4.20}]{Meereboer2025}\label{cor:herm2}
	Let $(\I_\bullet,\tau)$ be an irreducible Satake diagram. Then, the following three statements are equivalent:
	\begin{enumerate}[(i)]
		\item There exist countably infinite pairwise non-isomorphic integrable characters of $\uqbs$.
		\item There exist non-trivial integrable characters of $\uqbs$.
		\item $\uqbs$ is of Hermitian type.
	\end{enumerate}
\end{prop}
We understand the branching rules of integrable characters through the branching rules for the classical case. 

For each dominant weight $\lambda\in X^+$ let $L(\lambda)^1$ denote the highest weight module of weight $\lambda$ for the universal enveloping algebra $\mathbf{U}(\mathfrak{g})$ and let $\mathbf{U}(\mathfrak{g}^\Theta)\subset \mathbf{U}(\mathfrak{g})$ denote the universal enveloping algebra of $\mathfrak{g}^\Theta$.

\begin{thm}\cite[\mbox{Cor 4.19}]{Meereboer2025}\label{cor:hermetian} Let $(\boldsymbol{c},\boldsymbol{s})$ be a quasi-specializable parameter and $\lambda\in X^+$. Then,
	for each integrable character $\chi\in \widehat{\uqbs}$ there exists a character $\cl(\chi)\in \widehat{\mathbf{U}(\mathfrak{g}^\Theta)}$ such that the following two statements are equivalent
	\begin{enumerate}[(i)]
		\item  The $\uq$-module $L(\lambda)$ restricted to $\uqbs$ contains a module of type $\chi$.
		\item The $\mathbf{U}(\mathfrak{g})$-module $L(\lambda)^1$ restricted to ${\mathbf{U}(\mathfrak{g}^\Theta)}$ contains a module of type $\cl(\chi)$.
	\end{enumerate}
\end{thm}
%Recall that the multiplicity of any character $\chi\in \widehat{\uqbs}$ in $L(\lambda)|_{\uqbs}$ is at most one for all $\lambda\in X^+$, cf. \cite[\mbox{Thm 3.4}]{Meereboer2025}.
\begin{defi}[Spherical weights]\label{def:spherical}
	We set $$\mathcal X^+(\uq,\uqbs,\chi):=\{\lambda \in  X^+: [ L(\lambda)|_{\uqbs}:V_\chi]=1\},$$
	and refer to their elements as $\chi$-\textit{spherical weights}. 
\end{defi} 
Let $P(\Sigma)\subset X$ denote the weight lattice associated to the root system $\Sigma$. Set $P(2\Sigma)=\{2\lambda:\lambda\in P(\Sigma)\}$ and $P^+(2\Sigma)=P(2\Sigma)\cap X^+$, we recall from \cite[\mbox{Thm 3.4}]{Letzter2003} that 
\[\mathcal X^+(\uq,\uqbs,\epsilon)= P^+(2\Sigma).\]
\begin{defi}[Bottom of the well]\label{def:well}\belowdisplayskip=-12pt
	The \textit{bottom of the $\chi$-well} is the set
	\begin{align*}
		&\mathcal  B^+(\uq,\uqbs,\chi):=\\
		&\{\lambda \in  X^+: [ L(\lambda)|_{\uqbs}:V_\chi]=1\text{ and }\lambda -\mu\notin \mathcal X^+(\uq,\uqbs,\chi)\text{ for all } \mu\in P^+(2\Sigma)\setminus \{0\}\}
	\end{align*}\qedhere
\end{defi}
This terminology is adopted from \cite[\mbox{Def 9.4}]{Pezzini2023}. As an example, the bottom of the well for the trivial character consists of the zero weight. 
By Theorem \ref{cor:hermetian} and \cite[Prop 9.7]{Pezzini2023}, we conclude that:
\begin{thm}\label{thm:well}
	For each integrable character $\chi\in \widehat{\uqbs}$, the set of $\chi$-spherical weights can be written as
	$$ \mathcal X^+(\uq,\uqbs,\chi)= \mathcal X^+(\uq,\uqbs,\epsilon)+\mathcal B^+(\uq,\uqbs,\chi).$$
\end{thm}
\begin{notation}\label{not:l}
Let $(G,K)$ be a compact Hermitian symmetric pair, where $G$ is a connected simply connected complex semisimple Lie group. By \cite[\mbox{Thm 5.1.7.}]{Heckman1995} the group $\widehat{K}$
of unitary characters of $K$ is infinite cyclic, generated by, say, $\chi_1$. We denote its elements by $\chi_l=(\chi_1)^l$
for $l\in \Z$. We view each $\chi_l$
also as a character of the universal enveloping algebra $\mathbf{U}(\mathfrak{g}^\Theta)$. 
It follows from Theorem \ref{cor:hermetian}, it follows that specialization of characters is injective. Therefore, the integrable characters of $\uqbs$ can be parameterized such that 
\begin{equation}\label{eq:clasical}
	\cl(\chi^l)=\chi_l 
\end{equation}
holds {for each} $ l\in \Z.$
\end{notation}
Lemma \ref{cor:bottom} gives a description of the characters $\chi_l$ in terms of its embedding in irreducible highest weight $\uq$-modules. The possible action of $\uqbs$ on the corresponding one-dimensional module can be determined using \cite[Lemma 4.16]{Meereboer2025}.
Let $i\in \I_\circ$ and $(\uq_i,\uqbi)$ be the rank one QSP-coideal subalgebra introduced in Definition \ref{def:rankone}. 
We recall that $X_i=\mathrm{Hom}(Y_i,\Z)$ belongs to its root datum.
%\begin{rem}
%By a case-by-case analysis using Table \ref{table:1} we note that if the restricted root system $\Sigma$ is not reduced that there exists a unique $i\in \I_\circ $ such that $(\I_\bullet,\tau|_{\I_\bullet\cup\{i,\tau(i)\}})$ is of type $\mathsf{AIII_a}$.
%\end{rem} 
For the proof of Lemma \ref{cor:bottom}, we extend the Notation \ref{def:spherical} to the case $q=1$.
\begin{lem}\label{cor:bottom}
	Let $\uqbs$ be of Hermitian type and $i\in \I_{\scriptscriptstyle\otimes}$ as in Table \ref{table:1}. Then, for each integer $l\in \Z$ 
	\[\mathcal B^+(\uq,\uqbs,\chi_l)=\{\frac{|l|}{2}\sum\limits_{\substack{{\alpha\in \Sigma^+_{\ell}}}}\alpha +\mu^l_0\},\]
	for a $\mu^l_0\in X$ that depends on $l\in Z$, satisfies $\langle \alpha^\vee_j,\mu_0^l\rangle=0$ for all $j\in I\setminus I_{\scriptscriptstyle\otimes}$ and $\mu_0^l=0$ if $\Sigma$ is a root system of type $\mathsf{C}$.
%	for non-reduced root systems with $(\I_\bullet,\tau|_{\I_\bullet\cup\{i,\tau(i)\}})$ of type $\mathsf{AIII_a}$.
\end{lem}
\begin{proof}
By \cite[\mbox{ \S 3.1,  Thm 3.1, Prop 3.3}]{Ho02} and Theorem \ref{cor:hermetian}, it follows that
\[\mathcal X^+(\uq,\uqbs,\chi_l)\stackrel{Thm \ref{cor:hermetian}}{=}\mathcal X^+(U(\mathfrak{g}),U(\mathfrak{k}),\cl(\chi_l))\stackrel{\cite{Ho02} }{=}\frac{|l|}{2}\sum\limits_{\substack{{\alpha\in \Sigma^+_{\ell}}}}\alpha +\mu^l_0 +\mathcal X^+(\uq,\uqbs,\epsilon),\]
with $\mu_0^l\in \N[\omega _i,\omega_{\tau(i)}]$ and $\mu_0^l=0$ if $\mu_0^l=0$ if $\Sigma$ is a root system of type $\mathsf{C}$ (called Case $I$ in \cite[Thm 3.1]{Ho02}).
Thus, Theorem \ref{thm:well} implies that \[\bottom=\{\frac{|l|}{2}\sum\limits_{\substack{{\alpha\in \Sigma^+_{\ell}}}}\alpha +\mu^l_0\}.\qedhere\] 
\end{proof}
\begin{rem}
The result \cite[\mbox{ \S 3.1,  Prop 3.3}]{Ho02} is based on the seminal paper \cite{Schlichtkrull1984}, however the corresponding weights are made more explicit in \\\cite[\mbox{ \S 3.1,  Prop 3.3}]{Ho02}.
\end{rem}
\section{Orthogonal polynomials and spherical functions}\label{sec:OPSF}
\subsection{Macdonald-Koornwinder polynomials}\label{sec:macdo}
This section reviews the definition of Macdonald-Koornwinder polynomials. For a comprehensive treatment of these polynomials, see \cite{Macdonald2003} or \cite{Koornwinder91}. 

Let $\epsilon_1,\dots \epsilon_n$ be an orthonormal basis of a real Hilbert space $V$. The root system $\Phi$ of $\mathsf{BC_n}$ is the union of the vectors $\Phi=\Phi_s\cup \Phi_m\cup \Phi_{\ell}$, where
\begin{equation}\label{eq:bc}
\Phi_s=\{\pm\epsilon_1,\dots,\pm \epsilon_n\},\quad \Phi_m=\{\pm \epsilon _i\pm \epsilon_j\,:\, 1\leq i<j\leq n \},\quad \Phi_{\ell}=2\Phi_s.
\end{equation}
The Weyl group $W^\Phi$ is of type $\mathsf{C_n}$ and acts on the weight and root lattice $P(\Phi)$, $Q(\Phi)$ respectively.
\begin{notation}\label{not:A}
Let $A=\F[P(\Phi)]$ denote the group algebra of $P(\Phi)$ over $\F$. Again we note that $W^\Phi$ naturally acts on $A$. We write $A^{W^\Phi}$ for the space of invariants. 

For each weight $\lambda\in P(\Phi)$ we denote the corresponding basis element of $A$ by $e^\lambda$. Furthermore, denote $\F\llbracket P(\Phi)\rrbracket$ the space of formal power series in $ P(\Phi)$. 
\end{notation} 

\begin{enumerate}[(i)]
	\item There exists an involution $\overline{\,\cdot\,}:A\to A$ that satisfies
	\begin{equation}\label{eq:involution1}
		f=\sum_\lambda f_\lambda e^\lambda\mapsto \overline{f}=\sum_\lambda f_\lambda e^{-\lambda}.
	\end{equation}
	\item There exists a $\Q$-linear involution $^0:A\to A$ that satisfies
	\begin{equation}\label{eq:involution2}
		f=\sum_\lambda f_\lambda e^\lambda\mapsto f^0=\sum_\lambda \overline{f_\lambda} e^{\lambda}.
	\end{equation}
\end{enumerate} 
\begin{rem}
	This notation follows \cite[\mbox{\S 5.1}]{Macdonald2003}. The reader is cautioned to distinguish the involutions $\overline{\,\cdot\,}$ on $\F$ and $A$.
\end{rem}
We continue by introducing a family of weights. Let $a$ and $x$, $x_1,\dots x_m$ be indeterminates and define
\begin{equation}\label{eq:poch}
	(x;a)_{l}:=\prod_{i=0}^{l-i} (1-xa^i),\qquad(x;a)_{\infty}:=\prod_{i=0}^\infty (1-xa^i),\qquad (x_1,\dots ,x_m;a)_{\infty}:=\prod_{j=1}^m(x_j;a)_{\infty}.
\end{equation}
Each $k$-\textit{labeling} is determined by a tuple of scalars $k=(k_1,\dots ,k_5)\in \Q^5$. To each $k$-labeling, we associate the \textit{Macdonald-Koornwinder weight in base $q$ of type }$(\mathsf{C^\vee_n},\mathsf{C_n}) $
\begin{align}\label{eq:macweight}
	\triangledown_k=\triangledown:=\prod_{\alpha\in \Phi_s} \frac{(e^{2\alpha};q)_\infty}{(q^{k_1}e^\alpha,-q^{k_2}e^\alpha,q^{k_3+1/2}e^\alpha,-q^{k_4+1/2}e^\alpha;q)_\infty}\prod_{\alpha\in \Phi_m}\frac{(e^\alpha;q)_\infty}{(q^{k_5}e^\alpha;q)_\infty}
	\in \F\llbracket P(\Phi)\rrbracket.
\end{align}
\begin{rem}
Here, the weight (\ref{eq:macweight}) coincides with Koornwinder's weight \cite[\mbox{(5.1), (5.2)}]{Koornwinder91} by identifying $(a,b,c,d,t,q)=(q^{k_1},-q^{k_2},q^{k_3+1/2},-q^{k_4+1/2},q^{k_5},q)$.
\end{rem}
The constant term map is defined by
\begin{align}\label{eq:ct}
\ct:\F \llbracket P(\Phi)\rrbracket\to \F,\quad \sum_{\lambda\in P(\Phi)}a_\lambda e^\lambda\mapsto a_0,\quad \text{where}\quad \sum_{\lambda\in P(\Phi)}a_\lambda e^\lambda \in \F \llbracket P(\Phi)\rrbracket.
\end{align}
The weight gives rise to a bilinear form using the constant term.
The \textit{Macdonald-Koornwinder inner product} $\langle\,,\,\rangle: A\times A\to \F$ is given by
\begin{equation}
	\langle f,g\rangle_k=\ct(f \overline{g} \triangledown)/\ct(\triangledown),\qquad\qquad\text{for }f,g\in A.
\end{equation}
\begin{notation}
For $f\in A$, the notation $f=e^\lambda+l.o.t.$ indicates that $f=e^\lambda+\sum_{\mu< \lambda}a_\lambda e^\lambda$. 
\end{notation}
For each weight $\lambda\in X^+$ set $m_\lambda= \sum_{\mu \in W^\Phi\lambda }e^{\mu}$.
Let $\langle \,,\,\rangle$ denote the Macdonald-Koornwinder inner product. The bilinear form gives rise to an orthogonal basis of $A^{W^\Phi}$. 
\begin{thm}\cite[\mbox{Thm 5.5}]{koornwinder1992askey} \& \cite[\mbox{5.3.1, 5.3.4}]{Macdonald2003}\label{thm:charmacdo}\label{thm:charkoorn}
	There exists a unique basis $\{P_{\lambda}:\lambda\in P^+(\Phi)\}$ of $A^{W^\Phi}$ with the following two properties:
	\begin{enumerate}[(i)]
		\item For each $\mu,\lambda\in P^+(\Phi)$ with $\mu-\lambda\in Q^+(\Phi)\setminus \{0\}$ we have $\langle P_\lambda,m_\mu\rangle=0$.
		\item For each dominant weight $\lambda\in P^+(\Phi)$ it holds that $P_\lambda=e^\lambda+l.o.t.$ .
	\end{enumerate}
\end{thm}
%Whenever $\mu<\lambda$, the polynomials $P_\lambda$ and $P_\mu$ are orthogonal. Furthermore, orthogonality holds even when the weights $\lambda$ and $\mu$ are unequal and incomparable, cf. \cite[(5.3.4)]{Macdonald2003}.
We refer to the polynomials  $\{P_{\lambda}\,:\,\lambda\in P^+(\Phi)\}$ as \textit{Macdonald-Koornwinder polynomials}.
\subsection{Spherical functions of type $\chi$}\label{sec:spherical}
This section introduces the spherical functions associated with characters of quantum symmetric pairs and recalls their relation with Macdonald-Koornwinder polynomials.

Throughout this section, let $\uqbs,\uqds\subset \uq$ denote QSP-coideal subalgebras and let $\chi\in\widehat{\uqbs}$ and $\eta\in \widehat{\uqds}$ denote characters.
\begin{defi}[Spherical function]
	A $(\chi,\eta)$-\textit{spherical function} is a function $\varphi\in \qfa$ that transforms according to the following rule:
	$$\varphi(bxb')=\chi(b)\varphi(x)\eta(b')\qquad\text{for all }  b\in \uqbs,b'\in \uqds,\text{ and } x\in \uq.$$
The $(\epsilon,\epsilon)$-spherical functions are referred to as \textit{zonal spherical functions.}
\end{defi}
 This notion of zonal spherical functions is the analogue of the quantum zonal spherical functions of Letzter, cf. \cite[\S 4]{Letzter2004}, who however, works with left coideal subalgebras instead of right ones. The compatibility of the two approaches is established in Lemma \ref{lem:compar} below.

\begin{defi}\label{def:sph}
The $\F$-vector space of $(\chi,\eta)$-spherical functions for $(\uqds,\uqbs)$ is denoted by $^{\uqds}_{\quad\chi}{\qfa}_{\eta}^{\uqbs}$. 
\end{defi}

Let $\lambda\in X^+$ be a dominant weight and $\boldsymbol{a}\in (\F^\times )^{\I}$. For each weight $\mu=\lambda-\eta$, with $\eta=\sum_{i\in\I}t_i\alpha_i$ and $t_i\in \Z_{\geq0}$ set
$$\phi_{\boldsymbol{a}}^{L(\lambda)}(\eta):=\prod_{i\in\I}a_i^{t_i/2}.\qedhere$$
\begin{lem}\cite[\mbox{Lem 4.12}]{Meereboer2025}\label{lem:tech4} Let $\lambda\in X^+$ be a dominant weight and let $\boldsymbol{a}\in (\F^\times)^\I$, then the following two statements hold.
	\begin{enumerate}[(i)]
		\item The $\uq$-modules $L(\lambda)$ and $^{\Phi_{\boldsymbol{a}}}L(\lambda)$ are isomorphic, via $\Phi_{\boldsymbol{a}}:v_\mu\mapsto \phi^{L(\lambda)}_{\boldsymbol{a}}(\mu)v_\mu$ for weight vectors $v_\mu\in L(\lambda)_\mu$.
		\item The $\uq$-modules $L(\lambda)^\ast$ and $\big(^{\Phi_{\boldsymbol{a}}}L(\lambda)\big)^\ast$ are isomorphic, via $\big(\Phi^\ast_{\boldsymbol{a}}\big)^{-1}:f_\mu\mapsto (\phi^{L(\lambda)}_{\boldsymbol{a}}(\mu))^{-1}f_\mu$ for weight vectors $f_\mu\in L(\lambda)^\ast_\mu$
	\end{enumerate}
\end{lem}

Recall $\rho=\frac{1}{2}\sum_{\alpha \in \mathcal R^+}\alpha$ the half sum of the positive roots of $\mathfrak{g}$, cf. \cite[10.2]{Humphreys2010}. Furthermore, recall the isomorphism $\Phi_{\boldsymbol{a}}$ defined in (\ref{eq:aut}) and its action weight $\uq$-modules as defined in Lemma \ref{lem:tech4}. 

\begin{defi}[$\rho$-shift]\label{def:rho}
Let $\boldsymbol{\zeta}\in (\F^\times)^\I$ be such that $\zeta_i=q^{(2\rho,\alpha_i)}$ for each $i\in\I$ and define $\uqbs^\rho:= \Phi_{\boldsymbol{\zeta}}(\uqbs)$. For all $\lambda\in X^+$, we define the operators
\[\mathrm{Ad}(K_\rho):=\Phi_{\boldsymbol{\zeta}}:\uq\to \uq\qquad \text{and}\qquad  K_\rho:=\Phi_{\boldsymbol{\zeta}}:L(\lambda)\to L(\lambda).\qedhere\] 
%Here, the action on $\uq$ is via (\ref{eq:aut}) and the action on $L(\lambda)$ is defined using Lemma \ref{lem:tech4}. 
\end{defi}
We note that
\begin{equation}\label{eq:rhoshifted}\mathrm{Ad}(K_\rho )(x)K_\rho v=K_\rho x v,\qquad \qquad x\in \uq, v\in L(\lambda).
\end{equation}

In this paper we are interested in spherical functions which transform from the right through $\uqbs$ and from the left by $\uqbs^\rho$. This choice is motivated by the Weyl group invariance with respect to $(\uqbs^\rho,\uqbs)$ spherical functions, cf. \cite[\mbox{Thm 5.3}]{Letzter2003}. The significance of $\mathrm{Ad}(K_\rho)$ stems from the fact that $\mathrm{Ad}(K_\rho)$ squares to $S^{-2}$, cf. \cite[\mbox{\S 4.9 (1)}]{Jantzen1996}.

\begin{defi}\label{def:rho2}
Let $\chi\in\widehat{\uqbs}$ be a character and consider the associated one-dimensional $\uqbs$-module $V_\chi$. By \eqref{eq:rhoshifted}, $K_\rho V_\chi$ is a one-dimensional $\uqbs^\rho$-module. Denote $\chi_\rho:\uqbs^\rho\to \F $ the corresponding character of $\uqbs^\rho$. 
\end{defi}

Observe that the character $\chi_\rho$ satisfies \begin{equation}\label{lem:tech}
	\chi=\chi_{\rho}\circ \mathrm{Ad}(K_\rho).
\end{equation}
\begin{defi}[Restriction]
	Let $\Res_{\uq^0}: \qfa\to \mathrm{Hom}(\uq^0, \F)$ be the restriction to the subalgebra $\uq^0$. 
	Since $\uq^0$ is isomorphic to the group algebra of $Y$, we may identify $\mathrm{Hom}(\uq^0, \F)\cong \F[X]$. This defines a homomorphism 
	\[\Res:\qfa\to \F[X].\qedhere\]
\end{defi}
\begin{thm}\cite[\mbox{Thm 4.11}]{Meereboer2025}\label{thm:unique2}
	For each character $\chi$ of $\uqbs$ the multiplicity of $\chi$ in $L(\lambda)$ is at most one. If $\chi$ occurs in $L(\lambda)$, there exists a unique $(\chi_\rho,\chi)$-spherical function $\varphi\in L(\lambda)^\ast\otimes L(\lambda)$ for $(\uqbs^\rho,\uqbs)$ with $\Res(\varphi)=e^\lambda+l.o.t.$ .
\end{thm}
For an integrable character $\chi: \uqbs\to \F$, recall the bottom of the $\chi$-well $\bottom$ as introduced in Definition \ref{def:well}.
\begin{defi}[$\chi$-Spherical functions, Elementary \& Fundamental]\label{defi:elemn}
The elements of $^{\uqbs^\rho}_{\,\,\,\,\,\chi_\rho}{\qfa}_{\chi} ^{\uqbs}$ are called \textit{$\chi$-spherical functions.} A $\chi$-spherical function is called \textit{elementary} if it is contained in a single Peter-Weyl block. It is called \textit{fundamental} if it is both elementary and its corresponding weight lies at the bottom of the $\chi$-well.
\end{defi}
The QSP-coideal subalgebras from \cite{Letzter2003} are left coideals instead of right coideals subalgebras. The left and right QSP-coideal subalgebras are isomorphic via the Chevalley involution $\omega$, that sends the Hopf structure to the opposite Hopf structure, cf. \cite[\mbox{Sec 2}]{Letzter2000}. Let $\uqbl$ denote the left coideal subalgebra $\omega(\uqbs)$. The notion of spherical functions has a natural correspondence for the coideal subalgebras $\uqbs^\omega$.
In the following lemma, we show that these differences do not affect the restrictions of the zonal spherical functions. 

\begin{notation}
Let $w_0\in W$ denote the unique longest element of the Weyl group associated to the simply connected root datum $(Y,X,\langle\,,\,\rangle,\dots)$.
\end{notation}

Recall that there exists a diagram automorphism $\tau_0:\I\to \I$ with $-w_0\alpha_i=\alpha_{\tau_0(i)}$, for each $i\in\I$. For each dominant weight $\lambda \in  X^+$, the twisted module $^\omega L(\lambda)$ is isomorphic to $L(-w_0\lambda)=L(\tau_0(\lambda))$.

 \begin{notation}\label{not:family}
Let $\{\psi_\lambda:\lambda\in \Ppl\}$ denote the family of elementary zonal spherical functions for $\uqbs$. Its existence is assured by Theorem \ref{thm:unique2}.
\end{notation}

\begin{lem}\label{lem:compar}
	Let $\uqbs$ be of Hermitian type and let $\{\varphi_\lambda:\lambda\in \Ppl\}$ denote the family of elementary zonal spherical functions for $\uqbl$. Then $\Res(\varphi_\lambda)=\Res(\psi_\lambda)$ for all $\lambda\in \Ppl$.
\end{lem}
\begin{proof}
	Let $\lambda\in \Ppl$. Note that $\psi_{-w_0\lambda}\circ \omega\in L(\lambda)^\ast\otimes L(\lambda)$ is a spherical function for $\uqbs^\omega$. Thus, using Theorem \ref{thm:unique2} it follows that $\psi_{-w_0\lambda}\circ \omega=\varphi_\lambda$.
	The action of $\omega$ implies that 
	\begin{equation}\label{eq:compar}
	\Res(\varphi_\lambda)=\Res(\psi_{-w_0\lambda}\circ \omega)=\overline{\Res(\psi_{-w_0\lambda})}.
	\end{equation}
	As $W^\Sigma$ is of type $\mathsf{C}$, we have $-1\in W^\Sigma$, and note that the action of $-1\in W^\Sigma$ on polynomials coincides with that of $(\ref{eq:involution1})$. Thus using $\overline{\Res(\psi_{-w_0\lambda})}\in A^{W^\Sigma}$, $-1\in W^\Sigma$ and Theorem \ref{thm:unique2} we obtain $$e^\lambda+l.o.t.\stackrel{Thm \ref{thm:unique2}}{=}\Res(\varphi_\lambda)\stackrel{\eqref{eq:compar}}{=}\overline{\Res(\psi_{-w_0\lambda})}\stackrel{-1\in W^\Sigma}{=}\Res(\psi_{-w_0\lambda})=e^{-w_0\lambda}+l.o.t.\, .$$
	Therefore $e^\lambda=e^{-w_0\lambda}$, implying that $\lambda=-w_0\lambda$. As a consequence $\Res(\varphi_\lambda)=\Res(\psi_\lambda)$ for all $\lambda\in P^+(2\Sigma)$.
\end{proof}
Proposition \ref{prop:weylinv} allows us to identify zonal spherical functions with orthogonal polynomials.  
% where $\F[X]$ is the group algebra of the weight lattice $X$.
\begin{prop}\label{prop:weylinv}
	Restriction to $\uq^0$ defines an isomorphism of the space of zonal spherical functions
	$\spqfa$ with $\F[P(2\Sigma)]^{W^\Sigma}.$
\end{prop}
\begin{proof}This follows from Lemma \ref{lem:compar} and \cite[\mbox{Cor 5.4}]{Letzter2003}.
\end{proof}
Recall the restricted root system $\Sigma$ from \eqref{eq:restricted}. For reduced root systems $\Sigma$, zonal spherical functions arise as Macdonald polynomials. 
For $\alpha\in \Sigma$ set \[\mathrm{mult}(\alpha)= |\{\beta\in \mathcal R: \widetilde{\beta}=\alpha\}|,\]
 cf. \cite[\S6]{Letzter2004}.
The parameters of the resulting Macdonald polynomials dependent on the multiplicity function via
\begin{equation}\label{eq:multiplicities}
	k_l=\mathrm{mult}(\widetilde{\gamma_1})(\widetilde{\gamma_1},\widetilde{\gamma_1})(\gamma_1,\gamma_1)/2,\quad\text{and}\quad k_m=\mathrm{mult}(\widetilde{\gamma_2})(\widetilde{\gamma_2},\widetilde{\gamma_2})(\gamma,\gamma_2)/2
\end{equation}
where $\widetilde{\gamma_1}\in \Sigma_{\ell},\,\widetilde{\gamma_2}\in \Sigma\setminus\Sigma_{\ell}.$ 

Recall the distinguished indices $\I_{\scriptscriptstyle\otimes}$ from Table \ref{table:1}. Furthermore, recall the notation $\Sigma_s=\frac{1}{2}\Sigma_{\ell}$ and $\Sigma_m$ from Remark \ref{rem:ell}.
\begin{thm}\cite[\mbox{Thm 8.2}]{Letzter2004}\label{lem:case}
	Let $\Sigma$ be reduced and let $\uqb$ be of Hermitian type. Then, the family $\{\Res(\varphi_\lambda):P^+(2\Sigma)\}$ coincides with Macdonald-Koornwinder polynomials of type $(\mathsf{C^\vee},\mathsf{C}) $ in base $q_i^2$ with parameters $(k_l,k_l,0,0,k_m)$, with $i\in I_{\scriptscriptstyle\otimes}$.
\end{thm}
\begin{proof}
	By \cite[\mbox{Thm 8.2}, \S 9]{Letzter2004}, the family of polynomials $\{\Res(\varphi_\lambda):P^+(2\Sigma)\}$ coincides with Macdonald polynomials in base $q_i^2$ of type $\mathsf{C}$. 
	By specializing the parameters $(k_1,k_2,k_3,k_4,k_5)$ to $k=(k_1,k_1,0,0,k_5)$ in the Macdonald-Koornwinder weight of type $(\mathsf{C^\vee},\mathsf{C}) $ in base $q_i^2$, one obtains the Macdonald-Koornwinder weight
	\begin{align}\label{eq:macweight2}
	\triangledown_k=&\prod_{\alpha\in 2\Sigma_s} \frac{(e^{2\alpha};q_i^2)_\infty}{(q_i^{2k_1}e^\alpha,-q_i^{2k_1}e^\alpha,q_ie^\alpha,-q_ie^\alpha;q_i^2)_\infty}\prod_{\alpha\in 2\Sigma_m}\frac{(e^\alpha;q_i^2)_\infty}{(q_i^{2k_5}e^\alpha;q_i^2)_\infty}\nonumber\\
	&=\prod_{\alpha\in 2\Sigma_s} \frac{(e^{2\alpha};q_i^4)_\infty(q_i^2e^{2\alpha};q_i^4)_\infty}{(q_i^{4k_1}e^{2\alpha};q_i^4)_\infty(q_i^2e^{2\alpha};q_i^4)_\infty}\prod_{\alpha\in 2\Sigma_m}\frac{(e^\alpha;q_i^2)_\infty}{(q_i^{2k_5}e^\alpha;q_i^2)_\infty}\nonumber\\
	&\stackrel{2\Sigma_s=\Sigma_\ell}{=}\prod_{\alpha\in 2\Sigma_\ell} \frac{(e^{\alpha};q_i^4)_\infty}{(q_i^{4k_1}e^{\alpha};q_i^4)_\infty}\prod_{\alpha\in 2\Sigma_m}\frac{(e^\alpha;q_i^2)_\infty}{(q_i^{2k_5}e^\alpha;q_i^2)_\infty}.
	\end{align}
	Note that the weight \eqref{eq:macweight2} coincides with the Macdonald weight of type $\mathsf{C}$ in base $q_i^2$ for the parameters $(k_1,k_5)$, associated to the finite root system $2\Sigma$, cf. \cite[\mbox{(5.1.7), (5.1.28) (iii)}]{Macdonald2003}. According to \cite[\mbox{Thm 8.2}]{Letzter2004},
	$k_1=k_l$ and $k_5=k_m$, with $k_l,k_m$ defined in (\ref{eq:multiplicities}).
\end{proof}
\begin{rem}\label{rem:lit}
We recall that \cite{Letzter2004}, solely studies the zonal spherical functions for reduced root systems $\Sigma$. The remaining Hermitian types are $\mathsf{AIII_a}$, $\mathsf{DIII_b}$ and $\mathsf{EIII}$, cf. Table \ref{table:1}. For $\mathsf{DIII_b}$ a family of zonal functions is identified in \cite{Noumi1995}, but here $(\mathrm{Ad}(K_\rho) \uqb^\ast,\uqb)$ spherical functions are considered instead of $( \uqb^\rho,\uqb)$ spherical functions. Here $\ast$ corresponds to the compact real form, cf. \cite[\mbox{6.1.7 (57)}]{Klimyk1997}. Problematically, the coideal subalgebra for $\mathsf{DIII_b}$ related to the reflection equation of \cite[\mbox{Thm 3.1, (4.2)}]{Noumi1995} is not $\ast$-invariant. As far as we know, in the $\mathsf{EIII}$ case, the zonal spherical functions have not been identified with orthogonal polynomials. The zonal spherical functions for type $\mathsf{AIII_a}$ have been identified with Macdonald-Koornwinder polynomials, cf. \cite{Noumi1996}.
% In type $\mathsf{AIII}_a$ the zonal spherical functions are known to arise as Macdonald-Koornwinder polynomials, cf. 
However, \cite{Noumi1996} does not fall into the regime of Letzter's quantum symmetric pair coideal subalgebras used in \cite{Letzter2004}, as it is written in the setting of algebras associated to \textit{reflection equations}. These algebras can be related to Letzter's QSP-coideal subalgebras, cf. \cite[\mbox{Sec 6}]{Letzter1999}. 
\end{rem}
Recall the Satake diagram of type $\mathsf{AIII_a}$ from Table \ref{table:1}. Let $\sigma\in \Q$ and $\mathbf{B}_{\boldsymbol{c}_\sigma}$ be the QSP-coideal subalgebra of type $\mathsf{AIII_a}$, this algebra is associated to the reflection equation depending on a rational parameter $\sigma$ \cite[\mbox{(2.14)}]{Noumi1996}, cf. \cite[\mbox{Sec 6}]{Letzter1999}. 
In Appendix \ref{apen:B}, we elaborate on the value of the associated parameter $\boldsymbol{c}_\sigma\in (\F^\times)^{\I_\circ}$.

Recall the Satake diagram of type $\mathsf{AIV_m}$ from Example \ref{ex:1}. In the Theorem \ref{lem:case2}, we let $m\geq2$ denote the integer such that $\mathsf{AIV_m}$ is a rank one diagram contained in the Satake diagram of type $\mathsf{AIII_a}$. 
\begin{thm}\cite[\mbox{Thm 3.4}]{Noumi1996}\label{lem:case2}
	Restrictions of quantum zonal spherical functions for the Satake diagram of type $\mathsf{AIII_a}$
	for $(\mathbf{B}_{\boldsymbol{c}_\sigma}^\rho,\mathbf{B}_{\boldsymbol{c}_\sigma})$ are Macdonald-Koornwinder polynomials in base $q^2$ with parameters
	$(1/2,\sigma+1/2,m-1,-\sigma,1)$.
\end{thm}
\begin{proof}
	According to \cite[\mbox{Thm 3.4}]{Noumi1996}, using the restriction $\Res$, the $(\mathbf{B}_{\boldsymbol{c}_\sigma},\mathbf{B}_{\boldsymbol{c}_\sigma})$ zonal spherical functions arise as Macdonald-Koornwinder polynomials. Appendix \ref{apen:B} provides the translation to $(\mathbf{B}_{\boldsymbol{c}_\sigma}^\rho,\mathbf{B}_{\boldsymbol{c}_\sigma})$ zonal spherical functions and the values of the corresponding parameters.
\end{proof}
 \section{Constructing orthogonal polynomials}\label{sec:Const}
In this section we construct orthogonal polynomials associated with spherical functions of type $\chi$. For each integrable character $\chi$, we assign distinguished a family of orthogonal polynomials $\{P_\lambda^\chi:\lambda\in P(2\Sigma)^+\}\subset \F[P(2\Sigma)]^{W^\Sigma}$ and identify the corresponding orthogonality weight. The constructed polynomials encode the family of spherical functions. An analogous construction appears in \cite{Aldenhoven2016}, \cite{Aldenhoven16a} for higher dimensional representations of specific QSP-coideal subalgebras, or more generally in the classical setting in \cite{Pruijssen2017}.
\subsection{The structure of spherical functions}
Recall from Definition \ref{def:sph} that the space of spherical functions is denoted by $^{\uqds}_{\,\,\,\,\,\,\,\chi'}\qfa_\chi^{\uqbs}$. Moreover, recall $\Phi_\chi: \uqbs\to \Phi_\chi(\uqbs)$, the shift of base-point as introduced in Definition \ref{defi:shifted}. 
To understand the structure of the spherical functions, we need to understand their multiplicative structure. Their multiplicative structure is studied using the shift of base-point. 

\begin{rem}
In this section we extensively use Sweedler notation \eqref{eq:sweed} for the coproduct. We remind the reader that we write $\triangle (x)=x_{(1)}\otimes x_{(2)}$, for $ x\in \uq$.
\end{rem}
\begin{lem}\label{lem:multchar1}Let $\chi\in \widehat{\uqbs}$, $\chi'\in\widehat{ \uqds}$ be characters, and let $\eta\in \widehat{ \Psi_{\chi}(\uqbs)}$, $\eta'\in\widehat{ \Psi_{\chi'}(\uqds)}$ be characters. Then the multiplication map restricted to \\$ ^{\uqds}_{\,\,\,\,\,\,\,\chi'}\qfa_\chi^{\uqbs}\times{^{\Psi_{\chi'}(\uqds)}_{\,\,\,\quad\quad\eta'}\qfa_\eta^{\Psi_{\chi}(\uqbs)}}$ maps to $^{\uqds}_{\eta'\circ {\Psi_{\chi'}}}\qfa^{\uqbs}_{\eta\circ \Psi_{\chi}}$. 
\end{lem}
\begin{proof}{\allowdisplaybreaks
		Let $b\in \uqbs$, $b'\in \uqds$ and $x\in \uq$. Then for all $\varphi\in {^{\uqds}_{\,\,\,\,\,\,\,\chi'}\qfa_\chi^{\uqbs}}$ and $\psi \in{ ^{\Psi_{\chi'}(\uqds)}_{\,\,\,\quad\quad\eta'}\qfa_\eta^{\Psi_{\chi}(\uqbs)}}$ we have
		\begin{align*}
			(\varphi\cdot \psi) (b'xb)&=\varphi(b'_{(1)}x_{(1)}b_{(1)})\psi (b'_{(2)}x_{(2)}b_{(2)})\\
			&=\chi'(b'_{(1)})\varphi (x_{(1)})\chi(b_{(1)})\psi (b'_{(2)}x_{(2)}b_{(2)})\\
			&=\varphi (x_{(1)}) \psi (\chi'(b'_{(1)})b'_{(2)}x_{(2)}\chi(b_{(1)})b_{(2)})\\
			&=\varphi (x_{(1)})\psi (\Psi_{\chi'}(b')x_{(2)}\Psi_{\chi}(b))\\
			&=(\eta'\circ \Psi_{\chi'})(b')\varphi (x_{(1)}) \psi (x_{(2)})(\eta\circ \Psi_{\chi})(b)\\
			&=(\eta'\circ \Psi_{\chi'})(b')(\varphi \cdot \psi )(x)(\eta\circ \Psi_{\chi})(b).
	\end{align*}}
	Hence, $\varphi\cdot \psi $ is a spherical function of type $(\eta'\circ \Psi_{\chi'},\eta\circ \Psi_{\chi})$ for $(\uqds,\uqbs)$.
\end{proof}
\begin{notation}
In this section let $n$ denote the rank of the Satake diagram $(\I_\bullet,\tau)$, i.e. the number of $\tau$-orbits in $\I_\circ$, cf. Definition \ref{def:rank}.
\end{notation}
Recall the notion of an elementary spherical function from Definition \ref{defi:elemn} and the $\rho$-shifted QSP coideal subalgebra $\uqbs^\rho$ from Definition \ref{def:rho}. 
\begin{prop}\label{prop:strucepsilon}
	There exist elementary zonal spherical functions $\phi_1,\dots ,\phi_n\in{ {^{\uqbs^\rho}_{\quad\epsilon}\qfa_\epsilon^{\uqbs}}}$ such that
	$^{\uqbs^\rho}_{\quad\epsilon}\qfa_\epsilon^{\uqbs}=\F[\phi_1,\dots ,\phi_n]$.
\end{prop}
\begin{proof}
	By \cite[\mbox{Lem 4.1}]{Letzter2003}, we have
	$${^{\uqbs^\rho}_{\quad\epsilon}\qfa_\epsilon^{\uqbs}}=\bigoplus_{\lambda\in \Ppl} {^{\uqbs^\rho}_{\quad\epsilon}\qfa_\epsilon^{\uqbs}}\cap L(\lambda)^\ast\otimes L(\lambda).$$
	Moreover, Theorem \ref{thm:unique2} implies that ${^{\uqbs^\rho}_{\quad\epsilon}\qfa_\epsilon^{\uqbs}}\cap L(\lambda)^\ast\otimes L(\lambda)\cong \F$ for each $\lambda\in \Ppl$. 
	It follows from the proof of \cite[Lemma 4.2.1]{Watanabe2023} that the $\Ppl$ is a free finitely generated monoid with $n$ generators.
	Let $\mu_1,\dots \mu_n\in \Ppl$ be the generators of the monoid $\Ppl$. For each $1\leq i \leq n$ choose a nonzero zonal spherical function 
	\[\phi_i\in {^{\uqbs^\rho}_{\quad\epsilon}\qfa_\epsilon^{\uqbs}}\cap L(\mu_i)^\ast\otimes L(\mu_i)\]
	with $\Res(\phi_i)=e^{\mu_i}+l.o.t.$.
	As $\Psi_\epsilon=\mathrm{id}$, it follows by Lemma \ref{lem:multchar1} that \[\F[\phi_1,\dots ,\phi_n]\subset{ ^{\uqbs^\rho}_{\quad\epsilon}\qfa_\epsilon^{\uqbs}}.\]
	For each weight $\lambda\in \Ppl$ let $\varphi_\lambda\in L(\lambda)^\ast\otimes L(\lambda)$ be a nonzero $\epsilon$-spherical function with $\varphi_\lambda=e^\lambda+l.o.t.$, its existence is assured by Theorem \ref{thm:unique2}. By induction on the height of $\lambda$, we show $\varphi_\lambda\in \F[\phi_1,\dots ,\phi_n]$. Write $\lambda=\sum_{i=1}^n a_i\mu_i\in \Ppl$. Then $Q_\lambda:=\prod_{i=1}^n \varphi_{\mu_i}^{a_i}\in \F[\phi_1,\dots ,\phi_n]$ and Theorem \ref{thm:unique2} implies that $\Res(Q_\lambda)=e^\lambda+l.o.t.$. By the induction hypothesis we have $\varphi_\lambda-Q_\lambda \in \F[\phi_1,\dots ,\phi_n]$. Hence it follows that $\varphi_\lambda\in \F[\phi_1,\dots ,\phi_n]$, completing the proof.
\end{proof}
\begin{rem}
Proposition \ref{prop:strucepsilon} is a quantum analog of \cite[\S3]{Vretare1976}.
\end{rem}
Recall the shift of base-point from Definition \ref{def:shift}.
In Theorem \ref{thm:isomod} we need the following identity of the shift of base-point: 
\begin{equation}\label{eq:chishift}
	\epsilon\circ \Psi_\chi(b)=\epsilon(\chi(b_{(1)})b_{(1)})=\chi(b_{(1)})\epsilon(b_{(2)})=\chi(b_{(1)}\epsilon(b_{(2)}))=\chi(b),\qquad b\in \uqbs.
\end{equation}
Furthermore, recall for each character $\chi\in \widehat{\uqbs}$ the associated character $\chi_\rho\in \widehat{\uqbs^\rho}$ from Definition \ref{def:rho2}. 
\begin{thm}\label{thm:isomod}
	The left $^{\uqbs^\rho}_{\,\,\,\,\,\,\,\,\,\epsilon}{\qfa}_{\epsilon}^{\uqbs}$-module $^{\uqbs^\rho}_{\,\,\,\,\,\chi_\rho}{\qfa}_{\chi} ^{\uqbs}$ is free and of rank one.
\end{thm}
\begin{proof}
	By Lemma \ref{cor:bottom}, there exists a unique dominant weight $b_\chi\in X^+$ such that $$\bottom=\{b_\chi\}$$
	Theorem \ref{thm:unique2} and Theorem \ref{thm:well} assure that for each $\lambda\in \Ppl$ there exists a unique \[\varphi_{\lambda+b_\chi}\in {^{\uqbs^\rho}_{\,\,\,\,\,\chi_\rho}{\qfa}_{\chi} ^{\uqbs}}\cap L(\lambda+b_\chi)^\ast \otimes L(\lambda+b_\chi)\]
	with $\Res(\varphi_{\lambda+b_\chi})=e^{\lambda+b_\chi}+l.o.t.$. We write $\varphi_\chi=\varphi_{b_\chi}$. By Lemma \ref{lem:multchar1}, noting by \eqref{eq:chishift} that $\epsilon\circ\Psi_{\chi}=\chi$, it follows that $ ^{\uqbs^\rho}_{\,\,\,\,\,\,\,\,\,\epsilon}{\qfa}_{\epsilon}^{\uqbs}\cdot \varphi_{\chi}\subset {^{\uqbs^\rho}_{\,\,\,\,\,\chi_\rho}{\qfa}_{\chi} ^{\uqbs}}$. 
%	for elementary zonal spherical functions $\phi_1,\dots ,\phi_n$.  Using Theorem \ref{thm:unique2}, for each weight $\lambda\in \Ppl$ we let 
%	\[\varphi_{\lambda+b_\chi}\in L(\lambda+b_\chi)^\ast\otimes L(\lambda+b_\chi)\]
%	be a $\chi$-spherical function with $\Res(\varphi_{\lambda+b_\chi})=e^{\lambda+b_\chi}+l.o.t.$. 
	
	Again, by induction on the height of $\lambda$, we show $\varphi_{\lambda+b_\chi} \in \F[\phi_1,\dots ,\phi_n]\cdot \varphi_{\chi}$. 
	Consider the weights $\mu_1,\dots ,\mu_n\in  P^+(2\Sigma)$ and the corresponding zonal spherical function $\phi_{1},\dots ,\phi_{n}$, as in the proof of Proposition \ref{prop:strucepsilon}. We recall from Proposition \ref{prop:strucepsilon} that \[^{\uqbs^\rho}_{\,\,\,\,\,\,\,\,\,\epsilon}{\qfa}_{\epsilon}^{\uqbs}=\F[\phi_1,\dots ,\phi_n],\]
	Write $\lambda=\sum_{i=1}^n a_i\mu_i\in \Ppl$ and set \[Q_\lambda:=\bigg(\prod_{i=1}^n \phi_{i}^{a_i}\bigg)\cdot \varphi_\chi\in \F[\phi_1,\dots ,\phi_n]\cdot \varphi_\chi.\]
	Theorem \ref{thm:unique2} implies that $\Res(Q_\lambda)=e^{\lambda+b_\chi}+l.o.t.$. By the induction hypothesis we have $\varphi_{\lambda+b_\chi}-Q_\lambda \in \F[\phi_1,\dots ,\phi_n]\cdot \varphi_{\chi}$. Hence it follows that $\varphi_{\lambda+b_\chi}\in \F[\phi_1,\dots ,\phi_n]\cdot \varphi_\chi$. 
	We thus conclude that $^{\uqbs^\rho}_{\,\,\,\,\,\chi_\rho}{\qfa}_{\chi} ^{\uqbs}$ is generated over $^{\uqbs^\rho}_{\,\,\,\,\,\,\,\,\,\epsilon}{\qfa}_{\epsilon}^{\uqbs}$ by $\varphi_\chi$. 
	
	Next, we verify that the module is freely generated. Suppose that the finite sum 
	\[\sum_{m_1,\dots m_n\in \Z_{\geq0}} a_{m_1,\dots ,m_n}\phi_1^{m_1}\dots \phi_n^{m_n}\cdot \varphi_\chi=0,\]
	for scalars $a_{m_1,\dots ,m_n}\in \F$.
	Let $(m_1,\dots ,m_n)\in (\Z_{\geq0})^{n}$ be such that the weight $\sum_{i=1}^n m_i\mu_i$ is maximal and $a_{m_1,\dots ,m_n}\neq 0$. By considering a maximal term in the restriction
	\[\Res\bigg(\sum_{m_1,\dots m_n\in \Z_{\geq0}} a_{m_1,\dots ,m_n}\phi_1^{m_1}\dots \phi_n^{m_n}\cdot \varphi_\chi\bigg)=a_{m_1,\dots ,m_n}e^{m_1\mu_1+\dots+m_n\mu_n+b_\chi}+l.o.t.=0,\]
	it follows that $a_{m_1,\dots ,m_n}=0$. A contradiction. So \[\sum_{m_1,\dots m_n\in \Z_{\geq0}}  a_{m_1,\dots ,m_n}\phi_1^{m_1}\dots \phi_n^{m_n}=0.\]
	In conclusion, the left $^{\uqbs^\rho}_{\,\,\,\,\,\,\,\,\,\epsilon}{\qfa}_{\epsilon}^{\uqbs}$-module $^{\uqbs^\rho}_{\,\,\,\,\,\chi_\rho}{\qfa}_{\chi} ^{\uqbs}$ is free and of rank one.
\end{proof}
\begin{defi}\label{def:assosiate}
Let $\varphi_\chi\in{^{\uqbs^\rho}_{\,\,\,\,\,\chi_\rho}{\qfa}_{\chi}^{\uqbs}}$ be the distinguished generator of the left $^{\uqbs^\rho}_{\,\,\,\,\,\,\,\,\,\,\epsilon}{\qfa}_{\epsilon}^{\uqbs}$-module $^{\uqbs^\rho}_{\,\,\,\,\,\chi_\rho}{\qfa}_{\chi} ^{\uqbs}$ with $\Res(\varphi_\chi)=e^{b_\chi}+l.o.t.$ . Using Theorem \ref{thm:isomod}, we denote 
\[\widehat{\,\cdot\,}:{^{\uqbs^\rho}_{\,\,\,\,\,\chi_\rho}{\qfa}_{\chi}^{\uqbs}}\to ^{\uqbs^\rho}_{\,\,\,\,\,\,\,\,\,\,\epsilon}{\qfa}_{\epsilon}^{\uqbs},\qquad \varphi\mapsto\widehat{\varphi}\]
the unique isomorphism of $^{\uqbs^\rho}_{\,\,\,\,\,\,\,\,\,\,\epsilon}{\qfa}_{\epsilon}^{\uqbs}$-modules with
$
	\varphi=\widehat{\varphi}\cdot \varphi_\chi.
$
\end{defi}
%\begin{cor}\label{cor:assosiate}
%	The map $\varphi\mapsto \widehat{\varphi}$ is an isomorphism of left $^{\uqbs^\rho}_{\,\,\,\,\,\,\,\,\epsilon}{\qfa}_{\epsilon}^{\uqbs}$-modules.
%\end{cor}
To relate $\chi$-spherical functions to orthogonal polynomials, we need to relate them to zonal spherical functions.

Recall $\mathrm{Ad}(K_\rho):\uq\to\uq$ and the action of $K_\rho$ on $\uq$-modules as given in Definition \ref{def:rho}.
\begin{defi}\label{def:xi}
Define a linear map $\Xi: \qfa\times \qfa \to \qfa$ by 
\[(\varphi,\psi)\mapsto \varphi\cdot( \psi\circ\mathrm{Ad}(K_\rho)\circ S),\qquad \text{for}\qquad \varphi,\psi\in \qfa.\qedhere\] 
\end{defi}
For matrix coefficients $c_{f,v}^{L(\lambda)}$ and $c_{g,w}^{L(\mu)}$ the element $\Xi\big(c_{f,v}^{L(\lambda)},c_{g,w}^{L(\mu)}\big)\in \qfa $ equals
$c_{f,v}^{L(\lambda)}\cdot S(c_{gK_\rho,K_{-\rho}w}^{L(\mu)}).$
\begin{lem}\label{lem:chitozonal}
	The map $\Xi$ restricted to ${^{\uqbs^\rho}_{\,\,\,\,\,\chi_\rho}{\qfa}_{\chi}^{\uqbs}}\times  {^{\uqbs^\rho}_{\,\,\,\,\,\chi_\rho}{\qfa}_{\chi}^{\uqbs}}$ maps to \\$ {^{\uqbs^\rho}_{\,\,\,\,\,\,\,\,\,\epsilon}{\qfa}_{\epsilon}^{\uqbs}}$.
\end{lem}
\begin{proof}
	Let $b'\in \uqbs^\rho,b\in\uqbs$ and $x\in \uq$. Then we have{\allowdisplaybreaks
		\begin{align*}
			&\varphi\cdot(\psi\circ\mathrm{Ad}(K_\rho)\circ S)(bxb')
			\\
			&=\varphi(b'_{(1)}x_{(1)}b_{(1)}) \psi(K_{\rho}S(b_{(2)})S(x_{(2)})S(b_{(2)}')K_{-\rho} )\\
			&=\chi_{\rho}(b'_{(1)})\varphi(x_{(1)})\chi(b_{(1)}) \psi(K_{\rho}S(b_{(2)})S(x_{(2)})S(b_{(2)}')K_{-\rho} )\\
			&\stackrel{\eqref{lem:tech}}{=}\chi(K_{-\rho} b'_{(1)}K_{\rho})\varphi(x_{(1)})\chi_{\rho}(K_{\rho}b_{(1)}K_{-\rho}) \psi(K_{\rho}S(b_{(2)})S(x_{(2)})S(b_{(2)}') K_{-\rho})\\
			&=\varphi(x_{(1)}) \psi(K_{\rho}b_{(1)}K_{-\rho}K_{\rho}S(b_{(2)})S(x_{(2)})S(b'_{(2)}) K_{-\rho}K_{-\rho} b_{(1)}'K_{\rho})\\
			&=\varphi(x_{(1)}) \psi(K_{\rho}b_{(1)}S(b_{(2)})S(x_{(2)})K_{-2\rho} S^{-1}(b_{(2)}') b_{(1)}'K_{\rho})\\
			&=\epsilon(b)\varphi(x_{(1)}) \psi(K_{\rho}S(x_{(2)})K_{-\rho})\epsilon(b')\\
			&=\epsilon(b)\varphi\cdot (\psi\circ\mathrm{Ad}(K_\rho)\circ S)(x)\epsilon(b').
	\end{align*}}
	Here we use $S(x)K_{-2\rho}=K_{-2\rho}S^{-1}(x)$, which is equivalent to $S^2(x)=K_{-2\rho}xK_{2\rho}$, cf. \cite[\mbox{\S 4.9 (1)}]{Jantzen1996}.
	Hence, $\Xi(\varphi,\psi)$ is a zonal spherical function for $(\uqbs^\rho,\uqbs)$.
\end{proof}
\begin{rem}
	Aldenhoven considered a map analogous to $\Xi$ in \cite[\mbox{Thm 4.5.8}]{Aldenhoven16a}, however there a $\ast$-structure was used instead of the antipode $S$.
\end{rem}
To go further, we need to know more about $\Xi$. 
\begin{defi}[Haar integral]
Denote $h: \qfa\to L(0)^\ast \otimes L(0)\cong \F$ the projection map onto the trivial Peter-Weyl block, commonly known as the \textit{Haar integral}.
\end{defi}
 We can apply quantum Schur orthogonality to obtain orthogonality relations for the Haar integral. This idea is used in \cite{Koornwinder1993} and \cite{Noumi1992} to identify the Haar integral on the algebra $\,^{\uqbs^\rho}_{\quad\,\epsilon}\qfa_{\epsilon}^{\uqbs}$ with integration on a torus with respect to a weight. 
 Recall the basis of elementary zonal spherical functions $\{\psi_\lambda \,:\, \lambda \in P^+(2\Sigma)\}$ of $^{\uqbs^\rho}_{\quad\,\epsilon}\qfa_{\epsilon}^{\uqbs}$ with $\psi_\lambda\in L(\lambda)^\ast\otimes L(\lambda)$, as introduced in Notation \ref{not:family}.
 
\begin{lem}\label{lem:orthogo}
	The map $h\circ \Xi: \qfa\to \F$ satisfies
	\begin{enumerate}[(i)]
		\item $h\circ \Xi\bigg(L(\mu)^\ast \otimes L(\mu), L(\lambda)^\ast \otimes L(\lambda)\bigg)=0$ if $\lambda\neq \mu$.
		\item $h\circ \Xi\bigg(\psi_\lambda,\psi_\lambda\bigg)\neq 0$ for each $\lambda\in P^+(2\Sigma)$.
	\end{enumerate}
\end{lem}
\begin{proof}
	The first assertion follows from \cite[\mbox{11.2.2 Prop 15 (i)}]{Klimyk1997}. 
	
	Next, let $\uq_{\Q(q)}$ be the $\Q(q)$-subalgebra generated by $E_i,F_i$ and $K_h$, with $h\in Y$ and $i\in \I$. The $\uq_{\Q(q)}$-module $L(\lambda)_{\Q(q)}:=\uq_{\Q(q)}v_\lambda$ is a highest weight module of weight $\lambda$. By \cite[\mbox{Prop 4.1.3}]{Watanabe2023} and because the parameters are quasi-specializable, we may assume, without loss of generality, that the spherical vectors $v\in L(\lambda)$, $f\in L(\lambda)^\ast$ are normalized to lie in $L(\lambda)_{\Q(q)}$ and $L(\lambda)^\ast_{\Q(q)}$ respectively. Let $\{v_1,\dots v_m\}\subset L(\lambda)_{\Q(q)}$ be a basis consisting of weight vectors and  $\{f_1,\dots f_m\}$ be a dual basis.
	Remark that the diagonal operator $K_{-2\rho}=F=(F_{ij})$ is an intertwiner $L(\lambda)\to (L(\lambda)^\ast)^\ast$, cf. \cite[\mbox{5.3}]{Jantzen1996}. By \cite[\mbox{11.2.2 Prop 15 (2)}]{Klimyk1997} we have
	\begin{equation}\label{eq:invariant}
		h(c_{f_k,v_l} S(c_{f_{i},v_j}))=\delta_{k,j}\frac{F_{il}}{\mathrm{Tr}(F)}.
	\end{equation}
	Let $0\neq \psi_\lambda=c_{f,v}=\sum_{i,j=1} ^ma_{i,j}c_{f_i,v_j}\in L(\lambda)_{\Q(q)}^\ast\otimes L(\lambda)_{\Q(q)}$. Recall the ring $\A$ and the quotient homomorphism $\cl:\A\to \C$ as introduced in the beginning of Section \ref{sec:QG}. Choose $k\geq 0$ such that $(1-q)^ka_{i,j} \in \A$ for each $1\leq i ,j\leq m$ and there exist  $1\leq s,t \leq m$ with $\cl ((1-q)^ka_{s,t}) \neq 0$.
	Using (\ref{eq:invariant}) and $\cl(F_{i,l})=\delta_{i,l}$, it follows that
	\begin{align}\label{eq:trace}
		\cl\bigg(\mathrm{Tr}(F)h\circ \Xi((1-q)^k\psi_\lambda,(1-q)^k\psi_\lambda)\bigg)=\sum_{i,j=1} ^m \cl((1-q)^{k}a_{i,j})\cl((1-q)^ka_{j,i}).
	\end{align}
	For $1\leq i,j\leq m$ set $B_{i,j}=\cl((1-q)^{k}a_{i,j})$. As $B$ is nonzero, the trace of the positive matrix $BB^T$ is nonzero. Furthermore, $\mathrm{Tr}(B B^T)$ equals $(\ref{eq:trace})$.
	Thus, $h\circ \Xi\big(\psi_\lambda,\psi_\lambda\big)\neq 0$, completing the proof.
\end{proof}
\begin{rem}
	The strategy of showing that the form $h\circ \Xi$ is non-degenerate follows \cite[\mbox{Thm 3.7.}]{Dijkhuizen1994}. However, Dijkhuizen and Koornwinder's form is sesquilinear and complex-valued.
\end{rem}
\subsection{The construction}\label{sec:construction}
In this section we construct for each integrable character $\chi$ a distinguished family of orthogonal polynomials $\{P_\lambda^\chi:\lambda\in P(2\Sigma)^+\}\subset \F[P(2\Sigma)]^{W^\Sigma}$ and identify the corresponding orthogonality weight. Following Remark \ref{rem:lit}, we assume in this section that the Satake diagram is not of type $\mathsf{DIII_b}$ and $\mathsf{EIII}$. However, the constructions given in this section generalize to type $\mathsf{DIII_b}$ and $\mathsf{EIII}$ whenever their zonal spherical functions are identified with Koornwinder-Macdonald polynomials.

\begin{notation}\label{not:spherical}
Recall from Theorem \ref{thm:unique2} and Definition \ref{defi:elemn} that for each $\chi$-spherical weight $\lambda\in X^+$ there exists a unique $\chi$-spherical function $\varphi_\lambda\in L(\lambda)^\ast\otimes L(\lambda)$ with $\Res(\varphi_\lambda)=e^\lambda+l.o.t$. Definition \ref{def:assosiate} yields a family of zonal spherical functions $\{\widehat{\varphi}_\lambda:\lambda\in \Ppl\}$ satisfying
$$\varphi_{b_\chi+\lambda}=\widehat{\varphi}_\lambda\cdot \varphi_\chi\qquad\qquad \text{where }\lambda\in \Ppl\text{ and }b_\chi \in \mathcal B^+(\uq,\uqb,\chi).\qedhere$$
\end{notation}

\begin{defi}[Orthogonal polynomials]\label{defi:poly}For each weight $\lambda\in P^+(2\Sigma)$, set
\[
P^\chi_\lambda:= \Res (\widehat{\varphi}_{\lambda})=e^\lambda+l.o.t.
 .\]
Let $l\in \Z$ and consider the character $\chi^l$ from Notation \ref{not:l}, then we abbreviate
\[P_\lambda^{\chi^l}=P^l_\lambda.\qedhere\]
\end{defi}

\begin{rem}
For the rest of the paper, we set $A=\F[P(2\Sigma)]$. As in subsection \ref{sec:macdo} we write $A^{W^\Sigma}\subset A$ for the space of $W^\Sigma$-invariants.
\end{rem}

\begin{defi}[Bilinear forms]\label{def:bilinform}
	Define two bilinear forms on $^{\uqb^\rho}_{\,\,\,\,\epsilon}{\qfa}_{\epsilon}^{\uqb}$ by
	\begin{equation}\label{eq:bilin}
		\langle \varphi ,\psi \rangle= h\circ \Xi\big(\varphi ,\psi\big)\qquad\text{and}\qquad\langle \varphi ,\psi \rangle_\chi= h\circ \Xi\big(\varphi\cdot \varphi_\chi ,\psi\cdot \varphi_\chi\big),
	\end{equation}
	for $\varphi,\psi\in {^{\uqb^\rho}_{\,\,\,\,\,\,\epsilon}{\qfa}_{\epsilon}^{\uqb}}.$
	Extend these bilinear forms to $A^{W^\Sigma}$ via the isomorphism $\Res$ from Proposition \ref{prop:weylinv}, denoting them by the same symbols. 
\end{defi}
\begin{defi}[$\chi$-Shifted weight]\label{defi:shifted}
	Let $\triangledown$ be the Macdonald-Koornwinder orthogonality weight associated to the restriction of zonal spherical functions, cf. Theorem \ref{lem:case} \& \ref{lem:case2}. Let $\varphi_\chi$ be as in Definition \ref{def:assosiate}, and set $f_\chi=\Res(\varphi_\chi)$. 
Then the \textit{$\chi$-shifted weight} $\triangledown_\chi$ is defined by
\begin{equation}\label{eq:shiftedweifght}
	\triangledown_\chi=(f_\chi \overline{f_\chi})\triangledown.\qedhere
\end{equation}
\end{defi}
The $\chi$-shifted weight encodes the orthogonality for the polynomials \[\{P^\chi_\lambda:\lambda\in \Ppl\}\subset A^{W^\Sigma}.\]
For the following lemma, recall the $\{\psi_\lambda:\lambda\in\Ppl\}$ denote the basis of the zonal spherical functions, cf. Notation \ref{not:family}. Furthermore, recall from \eqref{eq:ct} the constant term map $\ct: A\to \F$.
\begin{lem}\label{lem:constant}
	For $f,g\in A^{W^\Sigma}$ we have $\langle f,g\rangle= \mathrm{ct}(f  \overline{g} \triangledown)/ \mathrm{ct}(\triangledown)$.
\end{lem}
\begin{proof}
	It follows from Lemma \ref{lem:orthogo}, that the bilinear form 
	\[\langle \,, \, \rangle:A^{W^\Sigma}\times A^{W^\Sigma}\to \F,\qquad \langle f,g\rangle = h\circ \Xi \big(\Res^{-1}(f),\Res^{-1}(g)\big),\qquad f,g\in A^{W^\Sigma} \]
	from Definition \ref{def:bilinform} 
	is non-degenerate and has the family $\{\Res(\psi_\lambda):\lambda\in P^+(2\Sigma)\}$ of restricted zonal spherical functions as an orthogonal basis.
	By Definition of $\Xi$ (Definition \ref{def:xi}) and Lemma \ref{lem:chitozonal}, we deduce that 
	\[\langle f,g\rangle =\langle 1, \overline{f}g\rangle,\qquad \text{for all }f,g\in A^{W^\Sigma}.\]
	Therefore, the bilinear form $\langle \,, \, \rangle:A^{W^\Sigma}\times A^{W^\Sigma}\to \F$ is uniquely determined by the linear form 
	\[l: A^{W^\Sigma}\to \F,\qquad f\mapsto \langle 1,f\rangle,\qquad \text{ for all }f\in A^{W^\Sigma}.\]
	We extend $l$ uniquely to a linear form $l: A\to \F$ by declaring that it is constant on ${W^\Sigma}$-orbits. Now, the element $\triangledown= \sum_{\mu\in P(2\Sigma)} e^{-\mu}l(e^\mu)$ satisfies
	\begin{align*}
	l(e^\mu)=\ct (e^\mu \triangledown'),\qquad \text{for all }\mu\in P(2\Sigma).
	\end{align*}
	Hence, we deduce that
	\[
	\langle f,g\rangle=\ct(f\overline{g}\triangledown'),\qquad \text{for all }f,g\in A^{W^\Sigma}.
	\]
	Recall form Definition \ref{defi:shifted}, that $\triangledown$ denotes the Macdonald weight associated to $\{\Res(\psi_\lambda):\lambda\in P^+(2\Sigma)\}$. 
	The linear maps $l':f\mapsto \ct(f\triangledown')$ and
	$l:f\mapsto\ct(f\triangledown)$ are both elements of
	$A^*=\mathrm{Hom}(A,\F)$.
	Let $\lambda\in P^+(2\Sigma)\setminus\set{0}$, then using Lemma \ref{lem:orthogo}, it follows that
	\[
	0 = h(\Xi(\psi_\lambda,\psi_0)) = 
	\overline{\psi_0(1)}\ct(\triangledown'\Res(\psi_\lambda))
	= \overline{\psi_0(1)}l'(\Res(\psi_\lambda)),
	\]
	as $\psi_0$ is a constant.
	Similarly, by Theorem \ref{lem:case} and Theorem \ref{lem:case2} $l(\Res(\psi_\lambda))=0$.
	We thus see that $l,l'$ have the same kernel of
	codimension 1: the span of $\psi_\lambda$ for $\lambda\ne0$.
	Consequently, $l,l'$ span the same 1-dimensional
	vector space, so they are linearly dependent.
	Since both are non-zero, they are proportional.
	Since $l,l'$ depend injectively on $\triangledown,\triangledown'$, we obtain that $\triangledown,\triangledown'$
	differ by a non-zero scalar.
	Evaluation of $l,l'$ at 1 gives that scalar.
\end{proof}
This translates into a description of the bilinear form $\langle \,,\,\rangle_\chi$ from Definition \ref{def:bilinform}.
\begin{prop}\label{cor:bilinform}
	For $f,g\in A^{W^\Sigma}$ we have $\langle f,g\rangle_\chi= \mathrm{ct}(f  \overline{g} \triangledown_\chi)/ \mathrm{ct}(\triangledown)$.
\end{prop}
\begin{proof}\belowdisplayskip=-12pt
	For $f,g\in A^{W^\Sigma}$ take, using Proposition \ref{prop:weylinv}, $\varphi,\psi\in {^{\uqb^\rho}_{\,\,\,\,\,\,\epsilon}{\qfa}_{\epsilon}^{\uqb}}$ with $\Res(\varphi)=f$ and $\Res(\psi)=g$. Recall that $\triangle\circ S=P\circ (S\otimes S)\circ \triangle$, where $P$ is the flip map, cf. \cite[3.8. (1)]{Jantzen1996}. Then we have{\allowdisplaybreaks
		\begin{align*}
			\langle f,g\rangle_\chi &=h\circ \Xi\big(\varphi\cdot \varphi_\chi ,\psi\cdot \varphi_\chi \big)\\
			&=h\bigg(\varphi\cdot \varphi_\chi \cdot\big((\psi\cdot \varphi_\chi )\circ \mathrm{Ad}(K_\rho)\circ S\big)\bigg)\\
			&=h\bigg(\varphi\cdot \varphi_\chi \cdot(\varphi_\chi \circ \mathrm{Ad}(K_\rho)\circ S)\cdot (\psi\circ \mathrm{Ad}(K_\rho)\circ S)\bigg)\\
			&\stackrel{Lem\, \ref{lem:chitozonal}}{=}h\circ \Xi\big(\varphi\cdot \underbrace{\varphi_\chi \cdot (\varphi_\chi \circ \mathrm{Ad}(K_\rho)\circ S)}_{\in{^{\uqb^\rho}_{\,\,\,\,\,\,\epsilon}{\qfa}_{\epsilon}^{\uqb}}} ,\psi \big)\\
			&\stackrel{Lem\, \ref{lem:constant}}{=} \mathrm{ct} \big(f(f_\chi \overline{f_\chi}) \overline{g} \triangledown\big)/\ct(\triangledown)\\
			&= \mathrm{ct} \big(f \overline{g} \triangledown_\chi\big)/\ct(\triangledown).
	\end{align*}}
\end{proof}
Using the quantum Schur orthogonality relations as in Lemma \ref{lem:orthogo} we obtain orthogonality relations for the polynomials $\{P^\chi_\lambda\,:\,\lambda\in P^+(2\Sigma)\}$ from Definition \ref{def:assosiate}.
\begin{cor}\label{cor:ortho}
	For each $\lambda,\mu \in \Ppl$ we have
	$\mathrm{ct}(P^\chi_\lambda \overline{P^\chi_\mu}\triangledown_\chi)=0$ if $\lambda\neq \mu$.
\end{cor}
\begin{proof}
Let $\varphi_{\lambda+b_\chi}$ and $\varphi_{\mu+b_\chi}$ be as in Notation \ref{not:spherical}. 
Proposition \ref{cor:bilinform} asserts that 
\[\mathrm{ct}(P^\chi_\lambda \overline{P^\chi_\mu}\triangledown_\chi)/ \mathrm{ct}(\triangledown)=h \circ \Xi(\varphi_{\lambda+b_\chi},\varphi_{\mu+b_\chi}),\]
that by Lemma \ref{lem:orthogo} is equal to zero if $\lambda\neq \mu$.
\end{proof}
\section{Fundamental spherical functions}\label{sec:fund}
Recall that the fundamental spherical functions are precisely those of the form $\varphi\in L(\lambda)^\ast\otimes L(\lambda)$ with $\lambda \in \bottom$. 
To identify the polynomials constructed in Section \ref{sec:construction}, it suffices, by Corollary \ref{cor:ortho}, to determine the fundamental spherical functions. 
This section determines all non-trivial fundamental spherical functions for standard QSP-coideal subalgebras of Hermitian type. We first calculate the non-trivial fundamental spherical functions in the rank one cases. The general result is obtained by reducing to the rank one case.
\subsection{Case of rank 1}\label{sec:rank1}
The rank-one Hermitian QSP-coideal subalgebras are of types $\mathsf{AIII_a}$ and $\mathsf{AIII_b}$, cf. Table \ref{table:1}. The rank one Satake diagram $\mathsf{AIII_a}$ is also denoted $\mathsf{AIV_n}$ with $n\geq 2$, and the rank one diagram of type $\mathsf{AIII_b}$ is also denoted by $\mathsf{AI_1}$. In the following two sections we determine, for each character, the nontrivial fundamental spherical functions. Recall, by Corollary \ref{cor:herm2} and Notation \ref{not:l} that the spherical functions are labeled by integers $l\in \Z$.
\subsubsection{Type $\mathsf{AI_1}$}\label{sec:AI}
The rank one QSP-coideal subalgebra $\uqbs$ of type $\mathsf{AI_1}$, is the subalgebra $\uqbs\subset \mathbf{U}_q(\mathfrak{sl}_2)$ generated by the element 
$B=F+cEK^{-1}+sK^{-1}$. 

The two dimensional module $L(\omega_1)$ has a basis $\{w_1,w_{2}\}$, the action of the generators of $ \mathbf{U}_q(\mathfrak{sl}_2)$ on this basis is given by
\begin{align}\label{eq:actionmod}
	K w_1&=qw_1\qquad\qquad &K w_{2}&=q^{-1}w_{2}\\
	E w_{1}&=0\qquad\qquad &  E w_{2}&= w_1\nonumber\\
	F{w_1}&= w_{2}\qquad\qquad & F w_{2}&=0.\nonumber
\end{align}
Recall the definition of the right $\uq$-module $L^r(\omega_1)= ^{\varrho}L(\omega_1)\cong L(\omega_1)^\ast$, as introduced in subsection \ref{subsec:qfa}.
Using \eqref{eq:actionmod}. The action of $\uqbs$ on the modules $L(\omega_1)$ and $L^r(\omega_1)$ is given by Table \ref{table:action}
{\begin{table}[H]
		\begin{center}
			\begin{tabular}{|c|c|c|}
				\hline
				& $w_1$ &$w_2$  \\
				\hline
				$B$	&$w_2+sq^{-1} w_1 $&$sqw_2+cqw_1$  \\
				\hline
				$\varrho(B)$	& $cqw_2+sq^{-1}w_1$  &$sqw_2+w_1$  \\
				\hline
			\end{tabular}
			\caption{The action of generators of $\uqbs$ and $\varrho(\uqbs)$ on $L(\omega_1)$}\label{table:action}
			\vspace*{-10mm}
		\end{center}
	\end{table}
	By Lemma \ref{lem:tech4} the spherical functions do not depend on the parameter $\boldsymbol{c}$. Without loss of generality, fix the parameter $c=q^{-1}$. For $s=0$ and $c=q^{-1}$, the vector $v=w_1+w_2$ is a spherical vector for $\uqb$. Lemma \ref{lem:tech4} implies that the vector $f=qK_{-\rho} v=w_1+qw_2\in L^r(\omega_1)$ is spherical of type $\chi^1$ for $\uqb^\rho$. Hence, using the identification $^{\varrho}L(\omega_1)\cong L(\omega_1)^\ast$, the matrix entry $c_{f,v}\in L(\omega_1)^\ast\otimes L(\omega_1)$ is a spherical function of type $\chi^1$ for $(\uqb^\rho,\uqb)$. 
	
	To proceed, we determine spherical vectors for  $\Psi_{\chi^l}(\uqb)$. Here, $\Psi_{\chi^l}$ denotes the shift of base-point, as introduced in Definition \ref{def:shift}.
	\begin{lem}\label{lem:form1}
		For each $l\geq 1$ the vector $v_l=q^{-l}w_1+w_2\in L(\omega_1)$ is $\chi^1$ spherical for
		$\Psi_{\chi^l}(\uqb)$ and the vector $f_l=q^{-l-1}w_1+w_2\in L^r(\omega_1)$ is $\chi^1$ spherical for 
		$\Psi_{\chi^l}(\uqb^\rho)$.
	\end{lem}
	\begin{proof}
		By Definition \ref{def:shift} $\Psi_{\chi^l}(\uqb)=\uqds$ is a QSP coideal subalgebra with parameters $d=q^{-1}$ and $t=[l]_q=\frac{q^l-q^{-l}}{q-q^{-1}}$. A direct computation, using Table \ref{table:action} and Lemma \ref{lem:tech4} shows that $v_l=q^{-l}w_1+w_2\in L(\omega_1)$ is $\chi^1$ spherical for 
		$\Psi_{\chi^l}(\uqb)$ and $f_l=q^{-l-1}w_1+w_2\in L^r(\omega_1)$ is $\chi^1$ spherical for 
		$\Psi_{\chi^l}(\uqb^\rho)$.
	\end{proof}
	Using the identification $L^r(\omega_1)\cong L(\omega_1)^\ast$ form subsection \ref{subsec:qfa}, we get by repeated application of Lemma \ref{lem:multchar1} that
	\begin{equation}\label{eq:phil}
	\varphi_l:=c_{f_1,v_1}\cdot \dots \cdot c_{f_l,v_l}\in (L(\omega_1)^{\otimes l})^\ast\otimes (L(\omega_1)^{\otimes l})
	\end{equation}
	is a spherical function of type $\chi^l$ for $\uqb$. According to \cite[\mbox{Prop 4.2}]{Klimyk1997} we have $\mathcal B_+(\uq,\uqb,\chi^l)=\{l\omega_1\}$. Since $L(l\omega_1)$ occurs with multiplicity one in 
	\[L(\omega_1)^{\otimes l}= L(l\omega_1)+\oplus_{n<l}m_n L(n),\qquad m_n\in \Z_{\geq0}\]
	and each weight $n\omega_1$ is not $\chi^l$-spherical for $n<l$, it follows that 
	$\varphi_l$ is the fundamental $\chi^l$-spherical function. 
	
	The Shapovalov form (\ref{eq:Shapovalof}) allows us to identify $L^r(\lambda)$ with $L(\lambda)^\ast$. Then, using Lemma \ref{lem:form1} and the notation from Notation \ref{not:A}, one computes that 
	\[\Res(c_{f_l,v_l})=\prod_{\alpha\in \Sigma^-}e^{\alpha/2}(1+q^{2l+1}e^{\alpha}).\]
	Thus, using \eqref{eq:phil}, the fundamental spherical function $\varphi_l$ satisfies
	\begin{equation}\label{eq:rank1r}
		\Res(\varphi_l)=\Res(c_{f_l,v_l}\cdot \dots \cdot c_{f_1,v_1})=-\prod_{\alpha\in \Sigma_{\ell}^-}e^{-l\alpha/2}(-q e^{\alpha};q^2)_l.
	\end{equation}
	Here, we remind the reader that $\Sigma_{\ell}\subset \Sigma$ denotes the set of long roots. Hence,
	\begin{equation}\label{eq:fundamentall}
	\Res(\varphi_l )\overline{\Res(\varphi_l)}=\prod_{\alpha\in \Sigma_{\ell}}(-q e^{\alpha};q^2)_l.
	\end{equation}
	To determine the fundamental spherical function associated to negative integers $-l$, note that  $\chi^{-l}(B)=[-l]_q$ and $\chi^{-l}=\chi^l\circ \Phi_{-}$. Here $\Phi_{-}$ corresponds to the automorphism (\ref{eq:aut}) with $\mathbf{a}=-1$. By \cite[\mbox{Thm 5.12}]{Meereboer2025}, $\Res(\varphi_{-l})=\Res(\varphi_l\circ \Phi_{-})$. Thus
	\begin{equation}\label{eq:rest}
		\Res(\varphi_{-l} )\overline{\Res(\varphi_{-l})}=\Res(\varphi_l)\overline{\Res( \varphi_l)}
		%=\prod_{\alpha\in \Sigma_{\ell}}(-q e^{\alpha};q^2)_l.
	\end{equation}
	From Lemma \ref{lem:tech4}, it follows that $\Res(\varphi_{-l} )\overline{\Res(\varphi_{-l})}$ also is given by $\eqref{eq:fundamentall}$.
	\subsubsection{Type $\mathsf{AIV_n}$, $n\geq 2$}\label{sec:aiv}
	The QSP-coideal subalgebra $\uqb$ of type $\mathsf{AIV_n}$ is the subalgebra $\uqb\subset \mathbf{U}_q(\mathfrak{sl}_{n+1})$ generated by
	$$K_1K_n^{-1},\quad B_1=F_1+c_1T_{w_\bullet}(E_{n})K_1^{-1},\quad B_n=F_n+c_nT_{w_\bullet}(E_{1})K_n^{-1},\quad \text{and}\quad \uq_\bullet,$$
	where $\I_\bullet=\{2,\dots n-1\}$ for $n>2$ and $\I_\bullet=\emptyset$ for $n=2$, and where $c_1,c_n\in \F^\times$. 
	
	We begin by determining spherical vectors in the vector representation $\mathbf{V}\cong L(\omega_1)$ of $ \mathbf{U}_q(\mathfrak{sl}_{n+1})$. Let $\{w_1,\dots ,w_{n+1}\}$ denote the canonical basis of $\mathbf{V}$. The generators $K_i$ act diagonally on these weight vectors. Explicitly their action is given by:
	\[K_i w_j=\begin{cases}qw_j&\text{if }i=j,\\
		q^{-1}w_j&\text{if }i=j-1\\
		w_j&\text{otherwise}.
	\end{cases}\]
	The nontrivial actions of the Chevalley generators $E_i,F_i$ are:
	\[E_iw_{i+1}=w_i,\qquad\qquad  F_iw_i=w_{i+1}\]
	and they annihilate all other basis vectors.
	
	\begin{lem}\label{lem:computation1}\label{lem:computation2}
		For $1\leq j\leq n$ and $i=1,n$ we have
		\begin{align*}&T_{w_\bullet}(E_{\tau(i)})(w_j)=\begin{cases}
				w_2&i=1,\, j=n+1;\\
				(-1)^{n}q^{-n+2}w_1& i=n,\,j=n;\\
				0&\text{else}.
			\end{cases}
		\end{align*}
		and
		\begin{align*}
			&\varrho(T_{w_\bullet}(E_{\tau(i)}))(w_j)=\begin{cases}
				w_{n+1}&i=1,\, j=2;\\
				(-1)^{n}q^{-n+2}w_{n}& i=n,\,j=1;\\
				0&\text{else}.
			\end{cases}
		\end{align*}
	\end{lem}
	\begin{proof}
		The first identity is shown in \cite[\mbox{Lem 4.2}]{Shen2023} for $n=2m+1$. An analogous proof works for even $n$. For the second identity we note that the Shapovalov form for $\mathbf{V}$ satisfies $(w_i,w_j)=\delta_{i,j}$. Using the $\varrho$-contravariance of the Shapovalov form and the fact that $\varrho(T_{w_\bullet}(E_{\tau(i)}))(w_j)$ is a weight vector, the conclusion follows. 
	\end{proof}
	\begin{lem}\label{lem:vector}
		The vector $v=c_1w_1-w_{n}\in \mathbf{V}$ is $\chi^1$-spherical for $\uqb$ and the vector $v^r=c_n^{-1}(-1)^{n}qw_{1}-w_{n+1}\in {^\varrho\mathbf{V}}$ is $\chi^1$ spherical for $\uqb^\rho$.
	\end{lem}
	\begin{proof}
		This is a direct computation using Lemma \ref{lem:computation2} and the description of the vector representation, as described in the beginning of this section.
	\end{proof}
	\begin{lem}\label{lem:form2}
		For each $l\geq 1$ the vector $v_l=q ^{-l}c_1w_1-w_{n}\in \mathbf{V}$ is $\chi^1$ spherical for 
		$\Psi_{\chi^l}(\uqb^\rho)$ and the vector $f_l=q^{-l}c_n^{-1}(-1)^nqw_{1}-w_n\in {^\varrho\mathbf{V}}$ is $\chi^1$ spherical for 
		$\Psi_{\chi^l}(\uqb^\rho)$.
	\end{lem}
	\begin{proof}
		Because $\chi^l(K_1K_{n}^{-1})=q^l$, Definition \ref{def:shift} implies that the QSP-coideal subalgebra $\Psi_{\chi^l}(\uqb)=\uqd$ has parameters  $d_1=c_1q^{-l}$ and $d_n=c_nq^l$. A direct computation, using Lemma \ref{lem:computation2} shows that $v_l=q ^{-l}c_1w_1-w_{n}$ is $\chi^1$ spherical for 
		$\Psi_{\chi^l}(\uqb^\rho)$ and $f_l=q^{-l}c_n^{-1}(-1)^nqw_{1}-w_n\in {^\varrho\mathbf{V}}$ is $\chi^1$ spherical for 
		$\Psi_{\chi^l}(\uqb^\rho)$.
	\end{proof}
	Using \ref{lem:form2} and the identification $^\varrho\mathbf{V}\cong \mathbf{V}^\ast$ from subsection \ref{subsec:qfa}, we get by repeated application of Lemma \ref{lem:multchar1} that
	$$\varphi_l:=c_{f_1,v_1}\cdot \dots \cdot c_{f_l,v_l}\in (L(\omega_1)^{\otimes l})^\ast\otimes (L(\omega_1)^{\otimes l})$$
	is a spherical function of type $\chi^l$ for $\uqb$. Lemma \ref{lem:multchar1} implies that $\mathcal B^+(\uq,\uqb,\chi^l)=\{l\omega_1\}$.
	By an argument analogous to that in subsection \ref{sec:AI}, $\varphi_l$ is the fundamental $\chi^l$-spherical function. 
	
	Recall that the Shapovalov form (\ref{eq:Shapovalof}) allows us to identify $L^r(\lambda)$ with $L(\lambda)^\ast$. Then, using Lemma \ref{lem:form1} and the notation from Notation \ref{not:A}, one computes that 
	\begin{equation}\label{eq:restriction}
		\Res(c_{f_l,v_l})=-\prod_{\alpha\in \Sigma_{\ell} }e^{\alpha/2}(1+c_1^{-1}c_n(-1)^{n} q e^{-\alpha}q^{2l}).
	\end{equation}
	We remind the reader that $\Sigma_{\ell}\subset \Sigma$ denotes the set of long roots. Thus, using (\ref{eq:restriction}), the spherical function $\varphi_l$ satisfies
	\begin{equation}\label{eq:rank1nr}
		\Res(\varphi_l)=-\prod_{\alpha\in \Sigma_{\ell} }e^{l\alpha/2}(c_1^{-1}c_n(-1)^{n} q e^{-\alpha};q^2)_l.
	\end{equation}
	Hence,
	\begin{equation}\label{eq:elemen}
		\Res(\varphi_l )\overline{\Res(\varphi_l)}=\prod_{\alpha\in \Sigma_{\ell}}(c_1^{-1}c_n(-1)^{n} q   e^{\alpha};q^2)_l.
	\end{equation}
	To determine the weight function for negative integers $-l$, note that $\chi^{-l}\circ\tau=\chi^{l}$, where $\chi^{-l}\circ \tau$ is a character of the QSP-coideal subalgebra $\tau(\mathbf{B}_{\boldsymbol{c}})$ and $\tau$ is part of the Satake diagram $(\I_\bullet,\tau)$. 
	Let $\vartheta_l\in \qfa$ be a fundamental spherical function of type $\chi^l$ for $\tau(\mathbf{B}_{\boldsymbol{c}})$. Then it holds that
	\begin{align}\label{eq:rest2}
		\Res(\varphi_{-l} )\overline{\Res(\varphi_{-l})}&=\Res( \vartheta_l\circ \tau)\overline{\Res( \vartheta_l\circ \tau)}\\
		&=\prod_{\alpha\in \Sigma_{\ell}}(c_1c_n^{-1}(-1)^{n} q e^{\alpha};q^2)_l.\nonumber
	\end{align}
	\subsection{A reduction to the case of rank 1}
	In this section we determine the fundamental spherical functions in general rank using a rank one reduction.
	
	\begin{rem}\label{rem:fix}
	Recall the marked indices $\I_{\scriptscriptstyle\otimes}\subset \I_\circ$, cf. Table \ref{table:1}. 
	For the remainder of this subsection, fix an $i\in \I_{\scriptscriptstyle\otimes}$. 
	\end{rem}
	Recall from Definition \ref{def:rankone} that the associated rank one QSP is $(\uqi,\uqbi)$. We perform a rank one reduction to this pair. 
	\begin{prop}
	The rank one subalgebra $\uqbi$ does not depend on the choice of marked index $i\in \I_{\scriptscriptstyle\otimes}$.
	\end{prop}
	\begin{proof}
	By Table \ref{table:1}, we see that there exist $i\in I_\circ$ such that $\I_{\scriptscriptstyle\otimes}=\{i,\tau(i)\}$. Since $\uqbi=\mathbf{B}_{\boldsymbol{c}_{\tau(i)},\boldsymbol{s}_{\tau(i)}}$, it follows that $\uqbi$ does not depend on the choice of marked index $i\in \I_{\scriptscriptstyle\otimes}$.
	\end{proof}
	Recall the homomorphism $\Res(Y_i): X\to X_i$ from Definition \ref{def:rankone} and the sets $\bottom$ and $\mathcal{B}_+(\uqi,\uqbi,\chi^l|_{\uqbi})$ from Definition \ref{def:well}.
	The bottoms of the wells of are related via the restriction $\Res(Y_i)$.
	\begin{lem}\label{lem:fun}
		The restriction of weights $\Res(Y_i)$ maps $ \mathcal{B}_+(\uq,\uqbs,\chi^l)$ to\\$ \mathcal{B}_+(\uqi,\uqbi,\chi^l|_{\uqbi})$.
	\end{lem}
	\begin{proof}
		According to Lemma \ref{cor:bottom} the bottom of the $\chi^l$-wells equal
		\begin{align*} \mathcal{B}_+(\uq,\uqbs,\chi^l)=\big\{\frac{|l|}{2}\sum_{\alpha\in \Sigma^+_{\ell}}\alpha+\mu_0^l \big\},\quad  \mathcal{B}_+(\uqi,\uqbi,\chi^l|_{\uqbi})=\big\{\frac{|l|}{2}\sum_{\alpha\in (\Sigma^+_i)_\ell}\alpha +\mu_0^l\big\},\end{align*}
		where $\mu_0^l\in X_i$. Because the the long roots in $\Sigma^+\setminus \Sigma^+_i$ restrict trivially to $Y_i$, the conclusion follows. 
	\end{proof}
	Recall that $\varphi_\chi$
	is the fundamental $\chi$-spherical function and the associated polynomial is $f_\chi=\Res(\varphi_\chi)$, see Definition \ref{defi:shifted}. 
	Lemma \ref{lem:fun} controls the restriction of the fundamental spherical function to the rank one QSP $(\uq_i,\uqbis)$. 
	In Lemma \ref{lem:decomposition} we use this information to decompose $f_\chi$. 
	\begin{defi}\label{def:palpha}
%	Let $i\in \I_{\scriptscriptstyle\otimes}$ be an distinguished index.
	For $\alpha\in \Sigma$ and $l\in \Z_{\geq0}$ set
	\begin{equation}\label{eq:casesl}
		p_\alpha=\begin{cases}
			e^{-l\alpha/2}(-q_\alpha e^{\alpha};q_i^2)_l,&\text{if }\Sigma\text{ is reduced};\\
			e^{l\alpha/2}(c_i^{-1}c_{\tau(i)}(-1)^{n} q_\alpha e^{-\alpha};q_i^2)_l&\text{if } \chi^l(K_iK_{\tau(i)}^{-1})=q^l;\\
			e^{l\alpha/2}(c_ic_{\tau(i)}^{-1}(-1)^{n} q_i e^{-\alpha};q_i^2)_l&\text{if }\chi^{-l}(K_iK_{\tau(i)}^{-1})=q^{-l}.
	\end{cases}\qedhere\end{equation}
	\end{defi}
	\begin{lem}\label{lem:decomposition}
		There exist $k_i\in \F[ e^\beta:\beta\in X\text{ and }\langle \widetilde{\alpha_i}^\vee,\beta\rangle =0]$ such that
		$f_\chi=p_{\widetilde{\alpha}_i} k_i,$ were $p_{\widetilde{\alpha}_i} $ is as defined in (\ref{eq:casesl}).
	\end{lem}
	\begin{proof}
		Let $i\in \I_{\scriptscriptstyle\otimes}$ be as in Remark \ref{rem:fix} and let $Y_i$ be as in Definition \ref{def:rankone}. 
		Consider the Zariski dense decomposition $Y\supset Y^\Theta\oplus Y^{-\Theta}$ into $\pm1$ eigenspaces of $\Theta$. 
		The transformation behavior on $Y^\Theta$ is fixed as $c_{f,v}$ is a $\chi$-spherical function. Thus, $\Res(c_{f,v})$ defined by its restriction to $Y^{-\Theta}$. 
		%Thus, it suffices to consider $h\in Y^{-\Theta}$.
		Consider the Zariski dense decomposition $(Y_i^{-\Theta})^\perp\oplus Y_i^{-\Theta}\subset Y^{-\Theta}$, where 
		\[(Y^{-\Theta}_i)^\perp=\{h\in Y^{-\Theta}: \langle h, \alpha_i\rangle=0\}.\]
		%Recall that $\Res(c_{f,v})\in \F[P(2\Sigma)]$
		Let $h\in (Y^{-\Theta}_i)^\perp$.
		%		and consider the decomposition $h=u+y$, with $\Theta(u)=u\in Y_i^\perp/2$  and $\Theta(y)=-y\in Y_i^\perp/2$. Let $\varphi=c_{f,v}$ denote the fundamental spherical function. 
		%		We recall that  $\Res(c_{f,v})\in \F[P(2\Sigma)]$ and that for each $\mu\in P(2\Sigma)\subset X$ we have $\frac{\mu-\theta(\mu)}{2}=\mu,$
		%		so that 
		%		\[\langle \frac{h-\Theta(h)}{2},\mu\rangle=\langle h,\mu\rangle. \]
		%		We thus note that $c_{f,K_uv}=\chi(h)c_{f,v}$, and $\chi(h)=1$ by \cite[\mbox{Lem 4.16.}]{Mee24}.
		Similar to the proof of \cite[\mbox{Thm 5.11}]{Meereboer2025}, $c_{f,K_hv}=c_{f,K_yv}$ is a spherical function for $\uqbi$. By Lemma \ref{lem:fun}, the restriction $c_{f,K_hv}|_{\uq_i}$ is a fundamental spherical function for $\uqbi$. According to \cite[\mbox{Thm 4.12}]{Meereboer2025}, there exists a scalar $\eta_h\in\F$ such that $c_{f,K_hv}|_{\uq_i}=\eta_h\varphi_i$, here $\varphi_i$ denotes the fundamental spherical function of type $\chi|_{\uqbi}$ for $\uqbi$. 
		Let $\pi_1,\pi_2$ and $\pi_3$ be the projection maps \begin{align*}
			\pi_1&:Y^\Theta\oplus(Y_i^{-\Theta})^\perp\oplus Y_i^{-\Theta}\to Y_i\\
			\pi_2&:Y^\Theta\oplus(Y_i^{-\Theta})^\perp\oplus Y_i^{-\Theta}\to (Y_i^{-\Theta})^\perp\\
				\pi_3&:Y^\Theta\oplus(Y_i^{-\Theta})^\perp\oplus Y_i^{-\Theta}\to Y^\Theta
		\end{align*}
		For each $h\in Y^\Theta\oplus(Y_i^{-\Theta})^\perp\oplus Y_i^{-\Theta}$ we have
		\begin{align*}
			c_{f,v}(K_h)=\chi(K_{\pi_3(h)})\varphi_i(K_{\pi_1(h)})\eta_{\pi_2(h)},\qquad \eta_{\pi_2(h)}\in \F.
		\end{align*}
		As $c_{f,v}\in\qfa$, the map $\eta:Y_i^\perp\to \F$ given by $h\mapsto \eta_h$ is represented by $h\mapsto k_i(h)$, where $k_i\in \F[ e^{\beta}:\beta \in X\,\text{and}\,\langle \widetilde{\alpha_i}^\vee,\beta\rangle =0]$. By (\ref{eq:rank1r}) (\ref{eq:rank1nr}) and (\ref{eq:rest2}) the map $\varphi_i:Y^{\Theta}\oplus Y_i^{-\Theta}\to \F$ defined by $y\mapsto \varphi(K_y)$ is represented by $p_{\widetilde{\alpha}_i}(y)$, with $p_{\widetilde{\alpha}_i}$ as defined in (\ref{eq:casesl}).
		As $(Y_i^{-\Theta})^\perp\oplus Y_i^{-\Theta}\subset Y^{-\Theta}$ is Zariski dense, it follows that $f_\chi=\Res(c_{f,v})$ decomposes as $f_\chi=p_{\widetilde{\alpha}_i}k_i$. 
	\end{proof}
	\begin{defi}
	Let $\F\big(P(2\Sigma)\big)$ denote the field of fractions of $\F[P(2\Sigma)]$. Note that the action of $W^\Sigma$ of $\F[P(2\Sigma)]$ induces an action of $W^\Sigma$ on $\F\big(P(2\Sigma)\big)$. Elements of $\F\big(P(2\Sigma)\big)^{W^\Sigma}$ are called $W^\Sigma$-\textit{symmetric.}
	\end{defi}
	For $\alpha\in \Sigma$, recall the polynomial $p_\alpha$ from Definition \ref{def:palpha}.
	\begin{thm}\label{thm:fundamental}
		Set 
		\begin{equation}\label{eq:decomppa}
			p=\prod_{\alpha\in \Sigma_{\ell}} p_\alpha
		\end{equation}
		then $f_\chi\overline{f_\chi}$ is a nonzero scalar multiple of $p$.
	\end{thm}
	\begin{proof}
		Let $i\in\I_{\scriptscriptstyle\otimes}$ be as in Remark \ref{rem:fix}. Set $k_i'=\prod\limits_{\substack{\alpha\in \Sigma_{\ell},\,   \alpha\neq \pm\widetilde{\alpha_i} }}p_\alpha$ and choose $p'\in \F\big(P(2\Sigma)\big)$ such that $f_\chi \overline{f_\chi}p'=p$. By Lemma \ref{lem:chitozonal} $f_\chi \overline{f_\chi}$ is the restriction of a spherical function. Thus, Proposition \ref{prop:weylinv} implies that $f_\chi \overline{f_\chi}$ is $W^\Sigma$-symmetric. Since both $f_\chi \overline{f_\chi}$ and $p$ are $W^\Sigma$-symmetric, $p'$ must also be $W^\Sigma$-symmetric. Let $p_{\widetilde{\alpha}_i}$ and $k_i$ be as in Lemma \ref{lem:decomposition}, then we have
		$$p_{\widetilde{\alpha}_i}\overline{p_{\widetilde{\alpha}_i}}k_i\overline{k_i}p'\stackrel{Lem\, \ref{lem:decomposition}}{=}f_\chi \overline{f_\chi}p'=p\stackrel{(\ref{eq:decomppa})}{=}p_{\widetilde{\alpha}_i}\overline{p_{\widetilde{\alpha}_i}}k_i',$$
		which shows that $k_i\overline{k_i}p'=k_i'$. Because $k_i'$ and $k_i\overline{k_i}$ are elements of \[\F[ e^{\beta}:\beta \in X\,\text{and}\,\langle \widetilde{\alpha_i}^\vee,\beta\rangle =0]\]
		we have $p'\in \F\big(e^{\beta}:\beta \in X\,\text{and}\,\langle  \widetilde{\alpha_i}^\vee,\beta\rangle =0\big)$. As $p'$ is $W^\Sigma$-symmetric and \[p'\in \F\big(e^{\beta}:\beta \in X\,\text{and}\,\langle \widetilde{\alpha_i}^\vee,\beta\rangle =0\big),\]
		$p'$ is a nonzero scalar.
	\end{proof}
	\begin{rem}
	Theorem \ref{thm:fundamental} allows us to normalize $f_\chi\overline{f_\chi}$ such that $f_\chi\overline{f_\chi}=p$. For the rest of the paper, we assume that this is the case.
	\end{rem}
	\begin{rem}
		Our rank one reduction method mirrors Letzter's approach for calculating radial parts of central elements, cf. \cite[\mbox{Lem 6.5, Thm 6.7}]{Letzter2004}.
	\end{rem}
	\begin{rem}
		The $q\to 1$ limit of the spherical function $f_\chi\overline{f_\chi}$ determined in Theorem \ref{thm:fundamental} yields the corresponding $q=1$ case, cf. \cite[\mbox{Thm 5.2.2.}]{Heckman1995}. Note that Heckman's convention uses the root system $R=2\Sigma$, which is always considered to be of type $\mathsf{BC}$, so the short roots in $R$ coincide with the long roots in $\Sigma$.
	\end{rem}
\section{Identification with Macdonald-Koornwinder polynomials}\label{sec:iden}
Let $\uqb$ be a standard QSP-coideal subalgebra of Hermitian type. From Proposition \ref{cor:herm2}, it follows that the integrable characters are labeled by integers. In Section \ref{sec:construction} we construct for each integer $l\in \Z$ a family of orthogonal polynomials $\{P_\lambda^{l}:\lambda\in P^+(2\Sigma)\}$, see Definition \ref{defi:poly}. The orthogonality weights are derived using Theorem \ref{thm:fundamental} and Corollary \ref{cor:ortho}. 
In this section we show that these polynomials coincide with Macdonald-Koornwinder polynomials.

We first assume that the restricted root system $\Sigma$ is reduced. 
Recall from Theorem \ref{lem:case} that the restriction of zonal spherical functions $\{P^0_\lambda:\lambda\in \Ppl\}$ are identified as Macdonald-Koornwinder polynomials of type $(\mathsf{C}^\vee,\mathsf{C}) $. Denote $k=(k_1,k_2,k_3,k_4,k_5)=(k_l,k_l,0,0,k_m)$ the corresponding parameter, cf. \cite[\mbox{Thm 8.2}]{Letzter2004} and Theorem \ref{lem:case}. 

\begin{thm}\label{thm:main}
	If $\Sigma$ is reduced and $\uqb$ is of Hermitian type, the family $\{P^l_\lambda:P^+(2\Sigma)\}$ are Macdonald-Koornwinder polynomials in base $q_i^2$ of type $(\mathsf{C}^\vee,\mathsf{C}) $ with parameters $(k_l,k_l,0,|l|,k_m)$, with $i\in \I_{\scriptscriptstyle\otimes}$.
\end{thm}
\begin{proof}
	%Recall from Theorem \ref{lem:case} that the polynomials $\{P^0_\lambda\,:\,\lambda\in P^+(2\Sigma)\}$ associated to the zonal spherical functions are Macdonald-Koornwinder polynomials with parameters $(k_l,k_l,0,|l|,k_m)$. 
	Denote $\triangledown_k$ the Macdonald-Koornwinder orthogonality weight corresponding to the parameters $k$, cf. \eqref{eq:macweight}, Theorem \ref{lem:case}.
	By Corollary \ref{cor:ortho} and Definition \ref{defi:shifted}, the polynomials $\{P^l_\lambda\,:\,\lambda\in P^+(2\Sigma)\}$ are orthogonal with respect to the shifted weight $\triangledown_{\chi^l}=f_{\chi^l}\overline{f_{\chi^l}}\triangledown_k$. Let $\triangledown_{k,l}$ denote the Macdonald-Koornwinder weight for the shifted parameters $(k_l,k_l,0,|l|,k_m)$.
	We remain to show $\triangledown_{k,l}=f_{\chi^l}\overline{f_{\chi^l}}\triangledown_k.$
	Using Theorem \ref{lem:case} and (\ref{eq:macweight}) we deduce that the shifted weight $\triangledown _{k,l}$ is given by
	\begin{align*}
		\triangledown _{k,l}&\stackrel{(\ref{eq:macweight})}{=}\triangledown_k\cdot 
		\prod_{\alpha\in \Sigma_{\ell}}(-q_i e^{\alpha};q_i^2)_l\\
		&\stackrel{Thm\, \ref{thm:fundamental}}{=}\triangledown_k\cdot f_\chi\overline{f_\chi}.
	\end{align*}
	Using Theorem \ref{thm:charmacdo} we conclude that the polynomials $\{P^l_\lambda\,:\,\lambda\in P^+(2\Sigma)\}$ are Macdonald-Koornwinder polynomials for the parameters $(k_l,k_l,0,|l|,k_m)$.
\end{proof}

Next, we consider a case where the restricted root system $\Sigma$ is not reduced.
Recall from Section \ref{sec:OPSF} the QSP-coideal subalgebra $\mathbf{B}_{\boldsymbol{c}_\sigma}$ of type $\mathsf{AIII}_a$ that depends on a non-negative rational parameter $\sigma$.
In this section we identify the polynomials $\{P^{l,\sigma}_\lambda:\lambda\in \Ppl\}$ from Definition \ref{defi:poly} with Macdonald-Koornwinder polynomials. Recall from Theorem \ref{lem:case2} that the zonal spherical functions 
are identified as Macdonald-Koornwinder polynomials with parameters $(1/2,\sigma+1/2,m-1,-\sigma,1).$
\begin{thm}\label{thm:main2}
	The family of orthogonal polynomials $\{P^{l,\sigma}_\lambda:P^+(2\Sigma)\}$ for $(\mathbf{B}^\rho_{\boldsymbol{c}_\sigma},\mathbf{B}_{\boldsymbol{c}_\sigma})$ are Macdonald-Koornwinder polynomials in base $q^2$ with parameters \begin{equation}
		(k_1,k_2,k_3,k_4,k_5)=
		\begin{cases}
			(1/2,\sigma+1/2-l,m-1,-\sigma,1).&\text{if }l<0\\
			(1/2,\sigma+1/2,m-1,-\sigma+l,1).&\text{if }l>0.
		\end{cases}
	\end{equation}
\end{thm}
\begin{proof}
	First suppose $l>0$. 
	%Recall from Theorem \ref{lem:case2} that the polynomials $\{P^{0,\sigma}_\lambda:P^+(2\Sigma)\}$ associated to the zonal spherical functions are Macdonald-Koornwinder polynomials with parameters $(1/2,\sigma+1/2,m-1,-\sigma,1)$.
	Denote $\triangledown_\sigma$ the Macdonald-Koornwinder weight corresponding to the parameters $(1/2,\sigma+1/2,m-1,-\sigma,1)$, cf. \eqref{eq:macweight}, Theorem \ref{lem:case2}. 
	By Corollary \ref{cor:ortho} and Definition \ref{defi:shifted}, the polynomials $\{P^{l,\sigma}_\lambda : \lambda \in {P^+(2\Sigma)}\}$ are orthogonal with respect to the shifted weight $\triangledown_{\sigma,\l}=f_{\chi^l}\overline{f_{\chi^l}}\triangledown_\sigma$. 
	Using Theorem \ref{lem:case}, equations (\ref{eq:macweight}),(\ref{eq:casesl}), and Lemma \ref{lem:relpar} we see that the shifted weight $\triangledown_{\sigma,l}$ from Definition \ref{defi:shifted}, for the parameters $(1/2,\sigma+1/2,m-1,-\sigma+l,1)$ equals
	\begin{align*}
		\triangledown_{\sigma,l}=\triangledown_\sigma\prod_{\alpha\in \Sigma_{\ell}} (-q^{-2\sigma+1}e^\alpha;q^2)_l\stackrel{Thm\, \ref{thm:fundamental}}{=}\triangledown_\sigma f_\chi\overline{f_\chi}.
	\end{align*}
	Using Theorem \ref{thm:charmacdo} we conclude that the polynomials $\{P_\lambda^{l,\sigma}\}_{\lambda\in P^+(2\Sigma)}$ are Macdonald-Koornwinder polynomials with parameters $(1/2,\sigma+1/2,m-1,-\sigma+l,1)$. A similar argument works in the case where $l<0$.
\end{proof}
\begin{rem}
For the remaining Hermitian symmetric spaces with non-reduced root system, $\mathsf{DIII_b}$ and $\mathsf{EIII}$, the theory of zonal spherical functions requires further development. A generalization of \cite{Letzter2004} to the non-reduced cases could simplify the identification of $\chi$-spherical functions.
\end{rem}
\section{Some applications}\label{sec:apl}
\subsection{$\imath$-bar involutions and orthogonal polynomials}
Recall the involutions (\ref{eq:involution1}) and (\ref{eq:involution2}) of $A= \F[P(2\Sigma)]$, and the bar involution (\ref{eq:barinv}) of $\F$.
In this section we show that the polynomials $\{P^\chi_\lambda\,:\,\lambda\in P^+(2\Sigma)\}$ constructed in Definition \ref{defi:poly} satisfy the relation
$$P^\chi_\lambda=\sum_{\mu\in P^+(2\Sigma)}a_{\lambda,\mu}^\chi m_\mu,\qquad\text{where}\qquad \overline{a_{\lambda,\mu}^\chi}=a_{\lambda,\mu}^\chi\text{ for all }\mu \in P^+(2\Sigma).$$
In case of a reduced root system $\Sigma$ this gives a new proof of the $q\to q^{-1}$ symmetry of the corresponding Macdonald polynomials, cf. \cite[\mbox{5.3.2}]{Macdonald2003}. This symmetry, on the level of quantum groups, can be interpreted using $\imath$-bar involutions and braid group operators, cf. \cite{Bao2018} \text{and} \cite{Wang2023}. 

To prove this invariance, we first need a preliminary result on the bar-involution on $\F (P(2\Sigma))$. Extend the involution $f\mapsto f^0$ on $A$ from \eqref{eq:involution2} in the natural way to $\F\big(P(2\Sigma)\big)$. Elements $f\in \F\big(P(2\Sigma)\big)$ with $f=f^0$ are called \textit{bar-invariant.}
\begin{lem}\label{lem:bar}
	If  $f\in \F\big(P(2\Sigma)\big)$ is bar-invariant, then $f^{-1}$ is bar-invariant,
\end{lem}
\begin{proof}
	Let $f\in \F\big(P(2\Sigma)\big)$. We see that $1=1^0=f^0 (f^{-1})^0=f  (f^{-1})^0,$ thus $(f^{-1})^0=f^{-1}.$
\end{proof}
\begin{defi}[Balanced and uniform parameters]
A parameter $\boldsymbol{c}\in (\F^\times)^{\I_\circ}$ is called \textit{balanced} if $c_i=c_{\tau(i)}$ for all $i\in\I_\circ$, and a parameter $\boldsymbol{c}\in (\F^\times)^{\I_\circ}$ is called \textit{uniform} if
$c_i=(-1)^{\alpha_i(2\rho_\bullet ^\vee)}q^{(\alpha_i,w_\bullet \alpha_{\tau(i)}+2\rho_\bullet)}\overline{c_{\tau(i)}}$ for each $i\in \I_\circ$.
\end{defi}
 These distinguished parameters, are used to obtain symmetries of scalar valued spherical functions, see \cite[\mbox{Prop 5.3}]{Meereboer2025}.
\begin{thm}
	Let $\boldsymbol{c}$ be a balanced parameter. Then the polynomials $\{P_\lambda^\chi:\lambda\in P^+(2\Sigma)\}$ satisfy
	$$P^\chi_\lambda=\sum_{\mu\in P^+(2\Sigma)}a_{\lambda,\mu}^\chi m_\mu,\qquad\text{where}\qquad \overline{a_{\lambda,\mu}^\chi}=a_{\lambda,\mu}^\chi\text{ for all }\mu,\lambda \in P^+(2\Sigma),$$
\end{thm}
\begin{proof}Let $\{\varphi_{\lambda+b_{\chi}}:\lambda\in \Ppl\}$ be the family of $\chi$ spherical functions as in Definition \ref{defi:poly}. Using \cite[\mbox{Lemma 4.12. (iv)}]{Meereboer2025}, it follows that the associated polynomials $P^\chi_\lambda$ do not depend on the choice of balanced parameter $\boldsymbol{c}.$ Thus, by Lemma \ref{lem:bar}
	it is sufficient to show $\Res(\varphi_{\lambda+b_\chi})^0=\Res(\varphi_{\lambda+b_\chi})$, for each spherical weight $\lambda$ in the case that $\boldsymbol{c}$ is a balanced uniform parameter.
	
	Let $c_{f,v}\in L({\lambda+b_\chi})^\ast\otimes L({\lambda+b_\chi})$ be a $\chi$ spherical function for $(\uqb,\uqb)$.
	Let $u\in Y$ be in the Zariski dense subset with decomposition $u=z+w$ where $\Theta(z)=-z$ and $\Theta(w)=w$. Recall the diagram automorphism $\tau_0:\I\to \I$ introduced in Section \ref{sec:spherical}, that satisfies $-w_0\alpha_i=\alpha_{\tau_0(i)}$ for each $i\in \I$. By a case-by-case check for the Satake diagrams of Table \ref{table:1}, one sees the involution $\tau_0$ satisfies $\tau_0\circ \tau=\mathrm{id}$. 
	
	Then \cite[\mbox{Prop 5.3}]{Meereboer2025} implies that there is a normalization of $c_{f,v}$ with
	\begin{align}\label{eq:proof4}
		c_{f,v}(K_z)&=c_{\overline{f},\overline{v}}(K_{w_\bullet w_0(z)})\\
		&=c_{\overline{f},\overline{v}}(K_{-w_\bullet \circ \tau_0 (z)})\nonumber\\
		&=c_{\overline{f},\overline{v}}(K_{\Theta(z)})\nonumber\\
		&=c_{\overline{f},\overline{v}}(K_{-z}).\nonumber
	\end{align} 
	Using (\ref{lem:tech}), we see that $c_{fK_{-\rho},v}$ is a spherical function of type $\chi$ for $(\uqb^\rho,\uqb)$. As $\Xi(c_{fK_{-\rho},v},c_{fK_{-\rho},v})\in{^{\uqb^\rho}_{\,\,\,\,\,\,\epsilon}\qfa_\epsilon^{\uqb}}$ and $K_w\in \uqb$, we have \begin{equation}\label{eq:winv}
		\Xi(c_{fK_{-\rho},v},c_{fK_{-\rho},v})(K_{z+w})=\Xi(c_{fK_{-\rho},v},c_{fK_{-\rho},v})(K_{z}),
	\end{equation}
	moreover
	\begin{align}\label{eq:kashi}
		\Xi(c_{fK_{-\rho},v},c_{fK_{-\rho},v})(K_{z})&\stackrel{Def \ref{def:xi}}{=}c_{f,v}(K_{-\rho} K_z)c_{f,v}(K_{-\rho} K_{-z})\\
		&\stackrel{\eqref{eq:proof4}}{=}c_{\overline{f},\overline{v}}( K_{\rho} K_{-z})c_{\overline{f},\overline{v}}( K_{\rho} K_{z}\nonumber)\\
		&\stackrel{\footnotesize\cite[\mbox{(7.3.4)}]{Kashiwara1993}}{=}\overline{\Xi(c_{fK_{-\rho},v},c_{fK_{-\rho},v})(K_{z})}.\nonumber
	\end{align}
	Recall that $-1\in W^\Sigma$ and its action on $A=\F[P(2\Sigma)]$ coincides with that of $(\ref{eq:involution1})$. Thus from (\ref{eq:kashi}), (\ref{eq:winv}) and the action of $-1\in W^\Sigma$ it follows that 
	\begin{align}\label{eq:proof}
	\Res\big(\Xi(c_{fK_{-\rho},v},c_{fK_{-\rho},v})\big)&\stackrel{(\ref{eq:kashi}), (\ref{eq:winv})}{=}\overline{\Res\big(\Xi(c_{fK_{-\rho},v},c_{fK_{-\rho},v})\big)}^0\nonumber\\
	&\stackrel{-1\in W^\Sigma}{=}\Res\big(\Xi(c_{fK_{-\rho},v},c_{fK_{-\rho},v})\big)^0.
	\end{align}
	Using Definition \ref{def:assosiate}, we write 
	\begin{equation}\label{eq:proof1}
		c_{fK_{-\rho},v}= \phi_\lambda\cdot \varphi_\chi
	\end{equation}
	 where $\phi_\lambda\in {^{\uqb^\rho}_{\,\,\,\,\,\,\epsilon}\qfa_\epsilon^{\uqb}}$. 
	 
	 Let $Q_\lambda=\Res(\phi_\lambda)$. Remark that $Q_\lambda$ is the restriction of a spherical function, which by Proposition \ref{prop:weylinv} is invariant under the action of the restricted Weyl group. Therefore, as $-1\in W^\Sigma$ we have 
	 \begin{equation}\label{eq:proof3}
	 \overline{Q_\lambda}=Q_\lambda.
	 \end{equation}
	 Again using that $-1\in W^\Sigma,$ we deduce that 
	\begin{align*}
		\Res\big(\Xi(c_{fK_{-\rho},v},c_{fK_{-\rho},v})\big)^0&\stackrel{\eqref{eq:proof}}{=}\Res\big(\Xi(c_{fK_{-\rho},v},c_{fK_{-\rho},v})\big)\\
		&\stackrel{\eqref{eq:proof1}}{=}f_\chi\overline{f_\chi} Q_\lambda \overline{Q_\lambda}\\
		&\stackrel{\eqref{eq:proof3}}{=}f_\chi\overline{f_\chi} Q_\lambda^2.
	\end{align*}
	For $\lambda=0$, we see that there exists a $s\in \F^\times$ with  $\big(sf_\chi\overline{f_\chi}\big)^0=\Res(sf_\chi\overline{f_\chi})$. Using Lemma \ref{lem:bar}, we deduce there is a normalization of $Q_\lambda$ with $Q_\lambda^2=\big(Q_\lambda^0\big)^2$. As a result $Q_\lambda=\eta Q_\lambda^0$, with $\eta=\pm 1$. Note that $Q_\lambda$ is a scalar multiple of $P_\lambda$, thus $(P^\chi_\lambda)^0=\xi P^\chi_\lambda$ for a scalar $\xi\in \F^\times$. Because $(e^\lambda)^0=e^\lambda$ and $(P^\chi_\lambda)^0=e^\lambda+l.o.t.$ it follows that $(P^\chi_\lambda)^0=P^\chi_\lambda$.
\end{proof}
\subsection{An identification of the Macdonald $q$-differential operator}
The goal of this section is to identify operators that arise from the center of $\uq$ with Macdonald $q$-difference operators. 

Denote $\uqcon$ the simply connected quantized enveloping algebra, cf. \cite[\mbox{3.2.10}]{Joseph1995}. 
The center of $\uqcon$ has a distinguished basis $\{c_\mu:\mu\in X^+\}$. Each central element $c_\mu$ acts on the simple module $L(\lambda)$ as scalar multiplication by $\sum_{w\in W}q^{(w\mu,\lambda+\rho)}$, cf. \cite[7.1.17]{Joseph1995}. 
\begin{defi}\label{defi:qdif}
Let $\mu\in X$. Then we define the operator 
\[D^\chi_\mu: \F[P(2\Sigma)]^{W^\Sigma}\to \F[P(2\Sigma)]^{W^\Sigma}\]
to be the unique operator making the following diagram commute:
\begin{equation}\label{diagram:com}
	\begin{tikzcd}
		& ^{\uqb^\rho}_{\,\chi_\rho}\qfa_{\chi}^{\uqb} \arrow[rr, "\psi\mapsto c_{\mu}\triangleright\psi"] &  & ^{\uqb^\rho}_{\,\chi_\rho}\qfa_{\chi}^{\uqb} \arrow[rd, "\varphi_\chi\cdot \varphi\mapsto \varphi"] &                                               \\
		^{\uqb^\rho}_{\,\,\,\,\,\epsilon}\qfa_{\epsilon}^{\uqb} \arrow[ru, "\varphi\mapsto \varphi_\chi\cdot \varphi"] \arrow[rd,rightarrow,"\Res","\thicksim" '] &                                                                                         &  &                                                                                            & ^{\uqb^\rho}_{\,\,\,\,\,\epsilon}\qfa_{\epsilon}^{\uqb} \\
		& {\F[P(2\Sigma)]^{W^\Sigma}} \arrow[rr,"D_{\mu}^\chi"]                                                  &  & {\F[P(2\Sigma)]^{W^\Sigma}} \arrow[ru,"\thicksim" ',leftarrow,"\Res"]                                                     &                                              
	\end{tikzcd}
\end{equation}
Here $c_\mu\triangleright  -:\qfa\to\qfa$ denotes the operator on $\qfa$ defined by
$$c_{f,v}\mapsto c_{f,c_\mu v},\qquad \text{where}\qquad c_{f,v}\in L^\ast(\lambda)\otimes L(\lambda),\, \lambda \in X^+.\qedhere$$
\end{defi}
\begin{lem}\label{lem:eigen}
	The orthogonal polynomials $\{P^\chi_\lambda\,:\, \lambda\in P^+(2\Sigma)\}$ from Definition \ref{defi:poly} form a joint eigenbasis for the operators $D^\chi_\mu$, with eigenvalues given by
	\[D^\chi_\mu P^\chi_\lambda=\bigg(\sum_{w\in W}q^{(w\mu,\lambda+b_\chi+\rho)}\bigg)P^\chi_\lambda\] for each $\lambda\in P^+(2\Sigma),$
	where $b_\chi\in \mathcal B^+(\uq,\uqb,\chi)$.
\end{lem}
\begin{proof}
Recall that the polynomials $\{P^\chi_\lambda\,:\, \lambda\in P^+(2\Sigma)\}$ are restrictions of elementary spherical functions. As $c_\mu$ acts on the simple module $L(\lambda)$ as scalar multiplication by $\sum_{w\in W}q^{(w\mu,\lambda+\rho)}$, the result is a direct consequence of Definition \ref{defi:qdif}.
\end{proof}
\begin{rem}\label{rem:reduced}
	\begin{enumerate}[(i)]
		\item For the remainder of this subsection, we assume that the restricted root system  $\Sigma $ is reduced and not of type $\mathsf{EVII}$. 
		\item For the rest of the section, fix a weight $\mu\in X$ such that $\widetilde{\mu}$ is a minuscule weight for $\Sigma$. The appendix of \cite{Letzter2004} gives explicit examples of such weights $\mu$. 
	\end{enumerate}

\end{rem}
Recall that the Macdonald polynomials $\{P^l_\lambda\,:\, \lambda\in P^+(2\Sigma)\}$ constructed in Section \ref{sec:Const} are eigenfunctions for the Macdonald $q$-difference operator 

\begin{equation}\label{eq:macdooperator}
\widetilde{D}^{k,l}_{2\tilde{\mu}} f= \sum_{w\in W^\Sigma}(\triangle_{k,l}^+)^{-1}w (T_{2\tilde{\mu}} \triangle_{k,l}^+ f),\qquad \text{where}\qquad f\in A^{W^\Sigma}= \F[P(2\Sigma)]^{W^\Sigma},
\end{equation}
where $T_\mu$ is the operator $e^\lambda\mapsto q^{(2\tilde{\mu},\lambda)}e^\lambda$ \text{for} $\lambda\in \Z\I$ and $\mu \in \Ppl$,
and where $\triangle_{k,l}^+$ is a certain weight function, cf. \cite[\mbox{(5.1.7)}]{Macdonald2003}. 
The eigenvalues corresponding to the Macdonald $q$-difference operator are:
\begin{equation}\label{eq:macdo}
	\widetilde{D}^{k,l}_{2\tilde{\mu}}P_\lambda^l= \sum_{w\in W^\Sigma}q^{(w\tilde{\mu}, \lambda+\rho_{k',l})}P_\lambda^l,\qquad \text{for each}\quad \,\lambda \in P^+(2\Sigma),\,\,l\in \Z ,
\end{equation}
cf. \cite[\mbox{(5.3.3)}]{Macdonald2003}. Here $\rho_{k',l}$ denotes the element defined in \cite[\mbox{(1.5.2)}]{Macdonald2003} that depend on the parameters from Theorem \ref{thm:main}. We remind the reader that these parameters depend on an integer $l\in \Z$.

We first match the eigenvalues of Lemma \ref{lem:eigen} and \eqref{eq:macdo} in the zonal spherical case $l=0$.
\begin{lem}\label{lem:natural}
Let $N$ denote index of $W^\Sigma$ in $W$ and let $\lambda\in P^+(2\Sigma)$. Then 
\[\sum_{w\in W}q^{(w\mu, \lambda+\rho)}=N\sum_{w\in W^\Sigma}q^{(w\tilde{\mu}, \lambda+\rho_{k',0})}.\]
\end{lem}
\begin{proof}
	In the zonal spherical case, the diagram (\ref{diagram:com}) consists of left multiplication by the central element $c_\mu$. Letzter has shown that the action of the operator $D^\epsilon_{\mu}$ coincides, up to a scalar multiple, with the action of the Macdonald $q$-difference operator $\widetilde{D}^{k,0}_{2\tilde{\mu}}$, cf. \cite[\mbox{Thm 8.2}]{Letzter2004}. 
	Thus, the eigenvalues associated to the operators  $\widetilde{D}^{k,0}_{2\tilde{\mu}}$ and $D^\epsilon_{\mu}$, relative to the basis $\{P^\epsilon_\lambda\,:\, \lambda\in P^+(2\Sigma)\}$ are proportional. By Lemma \ref{lem:eigen} and (\ref{eq:macdo}) this implies the existence of $N\in \F$ with
	\[\sum_{w\in W}q^{(w\mu, \lambda+\rho)}=N \sum_{w\in W^\Sigma}q^{(w\tilde{\mu}, \lambda+\rho_{k',0})}.\]
	Counting the terms in the sums $\sum_{w\in W^\Sigma}q^{(w\tilde{\mu}, \lambda+\rho_{k',0})}$ and $\sum_{w\in W}q^{(w\mu, \lambda+\rho)}$, we observe that $N$ coincides with the index of $W^\Sigma$ in $W$. 
\end{proof}
\begin{prop}\label{prop:coset}
There exist right coset representatives $w_1,\dots ,w_N$ of $W^\Sigma$ in $W$ such that $\widetilde{w_i\mu}=\widetilde{\mu}$ for each $1\leq i \leq N$. 
\end{prop}
\begin{proof}
Let $w_1,\dots ,w_N$ be right coset representatives of $W^\Sigma$. 
Without loss of generality, assume $\widetilde{w_i\mu}\geq0$ for all $1\leq i\leq N$, where these are considered as weights in $P^+(2\Sigma)$. We note that $w_i\mu=\mu-\beta$ for some $\beta\in \Z_{\geq0} \I$. Suppose that $\beta\neq 0$. As $\widetilde{w_i\mu}=\widetilde{\mu}- \widetilde{\beta}=\gamma\in P^+(2\Sigma)$ and because $\widetilde{\mu}$ is minuscule and $\widetilde{\mu}-\gamma\in Q^+(2\Sigma)$, this contradicts \cite[\mbox{(7.2)}]{Letzter2004}. Thus $\widetilde{\beta}=0$ and $\widetilde{w_i\mu}=\widetilde{\mu}$ for each $1\leq i \leq N$. 
\end{proof}

Recall that $\rho_{k',0}$ denotes the element defined in \cite[\mbox{(1.5.2)}]{Macdonald2003}, depending on the parameters $k$ of the Macdonald-Koornwinder polynomials associated to the zonal spherical functions.

\begin{lem}\label{lem:techeval}
Let $w_1,\dots ,w_N\in W$ be as in Proposition \ref{prop:coset}. Then, for each $w\in W^\Sigma$ and $1\leq i\leq N$ we have $(w\mu,\rho_{k',0})=(ww_i\mu,\rho).$
\end{lem}
\begin{proof}
	As $w\circ \Theta=\Theta\circ w$ for all $w\in W^\Sigma$, cf. \cite[\mbox{Rem 2.6}]{Dobson2019}, and  $\Theta(\lambda)=\lambda$,
	it follows from Lemma \ref{lem:natural} that
	\begin{align}\label{eq:zeros}
		N\sum_{w\in W^\Sigma}q^{(w\tilde{\mu},\rho_{k',0}+\lambda)}\stackrel{Lem \ref{lem:natural}}{=}\sum_{w\in W} q^{(w \mu, \lambda+\rho)}&=\sum_{w\in W^\Sigma} \sum_{i=1}^Nq^{(\widetilde{w w_i\mu}, \lambda)+(ww_i\mu,\rho)}\\&=\sum_{w\in W^\Sigma} \sum_{i=1}^Nq^{(w\tilde{\mu}, \lambda)+(ww_i\mu,\rho)}.\nonumber
	\end{align}
	Let $\lambda\in P^+(2\Sigma)$ be dominant such that all $(w \tilde{\mu},\lambda)$ with $w\in W^\Sigma$ are distinct (this choice can be made by avoiding a finite number of hyperplanes).
	By (\ref{eq:zeros}) the polynomial 
	$$F(x)=\sum _{w\in W^\Sigma}\big(Nq^{(w \tilde{\mu},\rho_{k',0})}-\sum_{i=1}^Nq^{(ww_i\mu,\rho)}\big)x^{( w\widetilde{\mu},\lambda)}$$
	vanishes for infinitely many values of $x$, namely $x=q^m$, $m\geq1$. Therefore $F(x)$ is identically zero, which implies that for $w\in W^\Sigma$ the equality \[Nq^{(w \tilde{\mu},\rho_{k,0})}=\sum_{i=1}^Nq^{(ww_i\mu,\rho)}\]
	holds and thus $(w\mu,\rho_{k,0})=(ww_i\mu,\rho)$ for all $1\leq i \leq N$.
\end{proof}
\begin{rem}
	The proof of Lemma \ref{lem:techeval} mirrors the approach of \cite[\mbox{4.5.8}]{Macdonald2003}.
\end{rem}
Recall by Remark \ref{rem:reduced} that we assume $\Sigma$ to be reduced. Thus, by Corollary \ref{cor:bottom}, the bottom of the $\chi^l$ well  from Definition \ref{def:well} is given by $\{\frac{|l|}{2}\sum_{\alpha\in \Sigma^+_{\ell}}\alpha\}$. 
\begin{thm}
	For each $l\in \Z$ the operator $D^{\chi^l}_\mu$ from Definition \ref{defi:qdif} is a scalar multiple of the Macdonald $q$-difference operator $\widetilde{D}^{k,l}_{2\tilde{\mu}}$ \eqref{eq:macdooperator}.
\end{thm}
\begin{proof}
	By \cite[\mbox{(1.5.2)}]{Macdonald2003}, we have $\rho_{k',l}=\rho_{k'}-\frac{1}{2}|l|\sum_{\alpha\in\Sigma^+_{\ell}}\alpha$. Let $w_1,\dots, w_N\in W$ be as in Proposition \ref{prop:coset}. 
	Note that $\Theta(\frac{1}{2}|l|\sum_{\alpha\in\Sigma^+_{\ell}}\alpha)=-\frac{1}{2}|l|\sum_{\alpha\in\Sigma^+_{\ell}}\alpha$.
	Thus, by Lemma \ref{lem:techeval} and Theorem \ref{thm:main}, we see that the eigenvalues of the operators $D^{\chi^l}_\mu$ and $\widetilde{D}^{k,l}_{2\mu}$ on the basis of Macdonald-Koornwinder polynomials $\{P^l_\lambda\,:\, \lambda\in P^+(2\Sigma)\}$ are related by
	{\allowdisplaybreaks
		\begin{align*}
			\sum_{w\in W}q^{(w\mu,\lambda+\rho+\frac{1}{2}|l|\sum_{\alpha\in \Sigma^+_{\ell}}\alpha)}&= \sum_{w\in W^\Sigma}\sum_{i=1}^Nq^{( ww_i\mu,\lambda+\rho+\frac{1}{2}|l|\sum_{\alpha\in \Sigma^+_{\ell}}\alpha)}\\
			&= \sum_{w\in W^\Sigma}\sum_{i=1}^Nq^{(\widetilde{ ww_i\mu},\lambda+\frac{1}{2}|l|\sum_{\alpha\in \Sigma^+_{\ell}}\alpha)+(ww_i\mu,\rho)}\\
			&\stackrel{Lem \ref{lem:techeval}}{=}\sum_{w\in W^\Sigma}\sum_{i=1}^N q^{(w\widetilde{ \mu},\lambda+\frac{1}{2}|l|\sum_{\alpha\in \Sigma^+_{\ell}}\alpha)+(w\mu,\rho_{k,0})}\\
			&=\sum_{w\in W^\Sigma}\sum_{i=1}^N q^{(w\widetilde{ \mu},\lambda+\frac{1}{2}|l|\sum_{\alpha\in \Sigma^+_{\ell}}\alpha)+(w\tilde{\mu},\rho_{k,0})}\\
			&=\sum_{w\in W^\Sigma}\sum_{i=1}^N q^{(w\widetilde{ \mu},\lambda+\rho_{k',0}+\frac{1}{2}|l|\sum_{\alpha\in \Sigma^+_{\ell}}\alpha)}\\
			&=N \sum_{w\in W^\Sigma}q^{(w\widetilde{ \mu},\lambda+\rho_{k',0}+\frac{1}{2}|l|\sum_{\alpha\in \Sigma^+_{\ell}}\alpha)}\\
			&=N \sum_{w\in W^\Sigma}q^{(w\widetilde{ \mu},\lambda+\rho_{k',l})}.
	\end{align*}}
	Therefore, using Lemma \ref{lem:eigen} and Theorem \ref{thm:main} the operators $\widetilde{D}^{k,l}_{2\mu}$ and $D^{\chi^l}_\mu$ are diagonal with respect to the distinguished basis $\{P_\lambda^l:\lambda\in P^+(2\Sigma)\}$ of Macdonald-Koornwinder polynomials, and their eigenvalues are proportional. This implies that $D^{\chi^l}_\mu$ is a scalar multiple of the Macdonald $q$-difference operator $\widetilde{D}^{k,l}_{2\mu}$.
\end{proof}
\begin{rem}
	For the Satake diagram of type $\mathsf{AIII_a}$ a similar argument can be made to identify an operator arising from (\ref{diagram:com}) with the $q$-differential operators from \cite{Koornwinder91}. 
\end{rem}
\subsection{Connection coefficients}
Recall $\mathcal X^+(\uq,\uqbs,\chi^{l})$, the set of $\chi^l$ spherical weights from Definition \ref{def:spherical}. To ease the notation, we for now assume $l\geq0$ (the results extend to $l\leq 0$). Let \[\{\varphi^{l,(\boldsymbol{c},\boldsymbol{s})}_{\lambda}:\lambda\in \mathcal X^+(\uq,\uqbs,\chi^{l})\}\qquad  \text{and}\qquad  \{P^{l,(\boldsymbol{c},\boldsymbol{s})}_{\lambda}:\lambda\in P^+(2\Sigma)\}\]
denote the family of $\chi^l$ spherical functions and the corresponding polynomials for $\uqbs$, cf. Theorem \ref{thm:unique2} and Definition \ref{defi:poly}.
In past sections we suppressed the parameter $(\boldsymbol{c},\boldsymbol{s})$ from the notations of these families, its explicit inclusion is now needed. Let 
\[\Psi_{\chi^l}: \uqbs\to \Psi_{\chi^l}(\uqbs)\]
denote the shift of base-point from Definition \ref{def:shift} and let $\Psi_{\chi^l}(\boldsymbol{c},\boldsymbol{s})$ denote the parameters of $\Psi_{\chi^l}(\uqbs)$. 

Analogous to the rank one case, cf. Section \ref{sec:rank1}, the fundamental spherical function $\varphi_{b^{l}}^{l,(\boldsymbol{c},\boldsymbol{s})}$ has a decomposition
\begin{equation}\label{eq:funddecomp}
	\varphi_{b_{l}}^{l,(\boldsymbol{c},\boldsymbol{s})}=c_{f_1,v_1}\cdots c_{f_l,v_l},\qquad \qquad l\in \Z,\,b_l\in \mathcal B^+(\uq,\uqbs,\chi^l)
\end{equation}
where each $c_{f_i,v_i}$ is the fundamental $\chi^{1}$ spherical function for $\rho_{\chi^i}(\uqbs)$. By Lemma \ref{lem:multchar1}, for any $m\geq 1$ and $\lambda\in\mathcal X^+(\uq,\uqbs,\chi^{l})$, we have
\begin{align}\label{eq:lspherical}
	\varphi^{l,(\boldsymbol{c},\boldsymbol{s})}_\lambda c_{f_{l+1},v_{l+1}}\cdots c_{f_{l+m},v_{l+m}}&\stackrel{Lem\, \ref{lem:multchar1}}{=}\sum_{\mu\in\mathcal X^+(\uq,\uqbs,\chi^{l+m\,l})}a^{l,\lambda}_\mu \varphi^{l+m,(\boldsymbol{c},\boldsymbol{s})}_\mu,\\
	c_{f_{1},v_{l}}\cdots c_{f_{m},v_{m}}	\varphi^{l,\rho_{\chi^m}(\boldsymbol{c},\boldsymbol{s})}_\lambda &\stackrel{Lem\, \ref{lem:multchar1}}{=}\sum_{\mu\in\mathcal X^+(\uq,\uqbs,\chi^{l+m\,l})}d^{l,\lambda}_\mu \varphi^{l+m,(\boldsymbol{c},\boldsymbol{s})}_\mu.\nonumber
	\end{align}
Here, the scalars $a^{l,\lambda}_\mu, d^{l,\lambda}_\mu$ are commonly known as \textit{connection coefficients.}

Let $b_m$ denote the unique weight contained in the bottom of the $\chi^{m}$-well, cf. Lemma \ref{cor:bottom}. To understand (\ref{eq:lspherical}), we need to understand the tensor product decomposition $L(\lambda)\otimes L(b_m)$. 
\begin{defi}Let $\Lambda$ be the set of weights occurring in $L(b_m)$. 
The weights $\eta$ satisfying $[L(\lambda)\otimes L(b_m):L(\eta)]\geq1$ are contained in $\lambda+\Lambda$, cf. \cite[24.4]{Humphreys2010}. 
Set $\varGamma_\lambda=\mathcal X^+(\uq,\uqbs,\chi^{l+m})\cap (\lambda+\Lambda)$
\end{defi}
Using (\ref{eq:funddecomp}) and (\ref{eq:lspherical}) we have
\begin{align*}
	P_\lambda^{l,(\boldsymbol{c},\boldsymbol{s})}&=\sum_{\mu\in \Gamma_\lambda-b^l}a^{l,\lambda}_\mu P^{l+m,(\boldsymbol{c},\boldsymbol{s})}_\mu\qquad\text{where}\qquad a^{l,\lambda}_\mu\in \F.\\
	P_\lambda^{l,\rho_{\chi^m}(\boldsymbol{c},\boldsymbol{s})}&=\sum_{\mu\in \Gamma_\lambda-b^l}d^{l,\lambda}_\mu P^{l+m,(\boldsymbol{c},\boldsymbol{s})}_\mu\qquad\text{where}\qquad d^{l,\lambda}_\mu\in \F.
\end{align*}
This relation will be made more explicit in some rank one cases.
\subsubsection{The case of $\mathsf{AI}$}
Since $b_{\chi^1}=b_{\chi^{-1}}=\omega$ and $\Lambda=\{\omega,-\omega\}$, we have $\varGamma_\lambda=\{\lambda,\lambda-\alpha\}$. Consequently, we obtain the following relation for Macdonald-Koornwinder polynomials, that are Askey-Wilson polynomials in rank one, cf. \cite[\mbox{(6.5.8)}]{Macdonald2003}:
\begin{align}
	P_\lambda^{l,(\boldsymbol{c},\boldsymbol{s})}&=P_\lambda^{l+1,(\boldsymbol{c},\boldsymbol{s})}+a^{l,\lambda}_{\lambda-\alpha}P^{l+1,(\boldsymbol{c},\boldsymbol{s})}_{\lambda-\alpha}\qquad \text{where}\qquad l\in \Z,\lambda\in P^+(2\Sigma),\\
	P_\lambda^{l,\rho_{\chi^1}(\boldsymbol{c},\boldsymbol{s})}&=P_\lambda^{l+1,(\boldsymbol{c},\boldsymbol{s})}+d^{l,\lambda}_{\lambda-\alpha}P^{l+1,(\boldsymbol{c},\boldsymbol{s})}_{\lambda-\alpha}\qquad \text{where}\qquad l\in \Z,\lambda\in P^+(2\Sigma).\nonumber
\end{align}
\subsubsection{The case of $\mathsf{AIV}_2$}
Since $b_{\chi^1}=\omega_1$ and $\Lambda=\{\omega_1,-\omega_2,-\omega_1+\omega_2\}$, we have $\varGamma_\lambda=\{\lambda,\lambda-(\omega_1+\omega_2)\}$. Consequently, we obtain the following relation for Askey-Wilson polynomials:
\begin{align}
	P_\lambda^{l,\sigma}&=P_\lambda^{l+1,\sigma}+a^{l,\lambda}_{\lambda-(\omega_1-\omega_2)}P^{l+1,\sigma}_{\lambda-(\omega_1+\omega_2)},\qquad \text{where}\qquad l>0,P^+(2\Sigma), \sigma\in \Q_{\geq0},\\
	P_\lambda^{l,\sigma+1}&=P_\lambda^{l+1,\sigma}+d^{l,\lambda}_{\lambda-(\omega_1-\omega_2)}P^{l+1,\sigma}_{\lambda-(\omega_1+\omega_2)},\qquad \text{where}\qquad l>0,P^+(2\Sigma), \sigma\in \Q_{\geq0}.\nonumber
\end{align}
A similar analysis for $b_{\chi^{-1}}$ results in the relation for Askey-Wilson polynomials:
\begin{align}
	P_\lambda^{l,\sigma}&=P_\lambda^{l-1,\sigma}+a^{l,\lambda}_{\lambda-((\omega_1-\omega_2))}P^{l-1,\sigma}_{\lambda-(\omega_1+\omega_2)},\qquad \text{where}\qquad l<0,P^+(2\Sigma),\sigma\in \Q_{\geq0},\\
	P_\lambda^{l,\sigma-1}&=P_\lambda^{l-1,\sigma}+d^{l,\lambda}_{\lambda-((\omega_1-\omega_2))}P^{l-1,\sigma}_{\lambda-(\omega_1+\omega_2)},\qquad \text{where}\qquad l<0,P^+(2\Sigma),\sigma\in \Q_{\geq0}.
\end{align}
For both cases, $\mathsf{AI}$ and $\mathsf{AIV_2}$, we obtain results concerning relations of connection coefficients for parameter shifts of Askey-Wilson polynomials. In particular, we find only two non-zero terms in the expansion. These relations also can be analytically deduced, as the connection coefficients themselves are explicitly known, cf. \cite[\mbox{(7.6.3), (7.6.9)}]{Gasper2004}. 
\begin{appendices}
	\section{Zonal spherical functions as Macdonald-Koornwinder}\label{apen:B}
	In  \cite{Noumi1996} the zonal spherical functions for the \textit{Noumi-Sugitani} coideal subalgebras of type $\mathsf{AIII_a}$ are identified with Macdonald-Koornwinder polynomials. These results are obtained in the context of \textit{Noumi-Sugitani} coideal subalgebras. In this appendix we elaborate on the relevant parameters. 
	
	Let $\sigma\in \Q_{\geq0}$ and consider the solution $J^\sigma$ to the reflection equation, as given in \cite[\mbox{2.14}]{Noumi1996}. According to \cite[\mbox{Cor 6.6}]{Letzter1999} there exists an associated QSP-coideal subalgebra $\mathbf{B}_{\boldsymbol{c}_\sigma}$ with corresponding parameters $\boldsymbol{c}_\sigma=(c_i)_{i\in I_\circ}$. The QSP-coideal subalgebra $\mathbf{B}_{\boldsymbol{c_\sigma}}^\omega=\omega(\mathbf{B}_{\boldsymbol{c_\sigma}})$ has the same zonal spherical functions as the coideal subalgebra from \cite{Noumi1996}. To match the convention with \cite{Noumi1996}, let $i\in\I_{\scriptscriptstyle\otimes}$ with	 $\chi_1(K_iK_{\tau(i)}^{-1})=q$.
	\begin{lem}\label{lem:relpar}
		For the parameters $\boldsymbol{c}_\sigma$ of the QSP-coideal subalgebra $\mathbf{B}_{\boldsymbol{c_\sigma}}^\omega$ it holds that $c_ic_{\tau(i)}^{-1}=(-1)^mq^{2\sigma}$.
	\end{lem}
	\begin{proof}
		Using \cite[\mbox{Cor 6.6}]{Letzter1999} the elements
		$$L_{1,2}^+-q^{\sigma}S(L^-_{2,\tau(i)}),\quad \text{and}\quad L_{\tau(i)-1,\tau(i)}^+-q^{-\sigma}S(L^-_{l,\tau(i)-1})$$
		are generators of $\uqb^\omega$, here $L_{i,j}^{\pm}$ denote the $L$-operators, cf. \cite[\mbox{8.5.2}]{Klimyk1997}. By \cite[\mbox{6.2.3 (65), 8.5.2}]{Klimyk1997}, it follows that $c_ic_{\tau(i)}^{-1}=(-1)^{m}q^{2\sigma}$.
	\end{proof}
	\begin{prop}\label{cor:askwil}
		The restriction of $(\mathbf{B}^\rho_{\boldsymbol{c_\sigma}},\mathbf{B}_{\boldsymbol{c_\sigma}})$ zonal spherical functions are Macdonald-Koornwinder polynomials with parameters \[(k_1,k_2,k_3,k_4,k_5)=(1/2,\sigma+1/2,m-1,-\sigma,1).\]
	\end{prop}
	\begin{proof}
		Let $\{\vartheta_\lambda:\lambda\in \Ppl\}$ denote the family of zonal spherical functions for $(\mathbf{B}_{\boldsymbol{c_\sigma}},\mathbf{B}_{\boldsymbol{c_\sigma}})$, cf. Theorem \ref{thm:unique2}. According to \cite[\mbox{Thm 3.4}]{Noumi1996} the restriction of $(\mathbf{B}_{\boldsymbol{c_\sigma}},\mathbf{B}_{\boldsymbol{c_\sigma}})$ 
		zonal spherical functions equal Macdonald-Koornwinder polynomials in base $q^2$ with parameters $(1/2,\sigma+1/2,m-1,-\sigma,1).$
		Using Proposition \ref{prop:strucepsilon} we can express $\Res(\vartheta_\lambda)$ as $ P_\lambda(\Res(\phi_1),\dots ,\Res(\phi_k))$ for elementary zonal spherical functions $\phi_1,\dots ,\phi_{k}$ and a polynomial $P_\lambda$. Since $\{\vartheta_\lambda\triangleleft K_{-\rho}:\lambda\in \Ppl\}$ is a scalar multiple of the family of zonal spherical functions for $(\mathbf{B}^\rho_{\boldsymbol{c_\sigma}},\mathbf{B}_{\boldsymbol{c_\sigma}})$, we see that $\Res(\vartheta_\lambda\triangleleft K_{-\rho})= P_\lambda(\phi_1\triangleleft K_{-\rho},\dots,\phi_k\triangleleft K_{-\rho})$, up to a scalar multiple. By setting  $\Res (\phi_i\triangleleft K_{-\rho}):=\phi_i'$ for $1\leq i\leq k$, we observe that $\Res(\vartheta_\lambda\triangleleft K_{-\rho})=P_\lambda(\phi_1',\dots ,\phi_k')$. We conclude that the family $\{\Res(\vartheta_\lambda\triangleleft K_{-\rho}):\lambda\in \Ppl\}$ consists of Macdonald-Koornwinder polynomials with the desired parameters.
	\end{proof}
\end{appendices}
%\begin{center}
%	$\textsc{Data Availability}$
%\end{center}
%Data sharing not applicable to this article as no data sets were generated or analysed during the current study.
%\setlength{\bibsep}{0.0pt}
\bibliographystyle{alpha}
\bibliography{Reff}
{\addcontentsline{toc}{section}{References}}
\end{document}